\documentclass[12pt]{iopart}
\usepackage{amsthm}
\usepackage{graphicx}
\usepackage{subfig}
\usepackage{amssymb}
\usepackage[linesnumbered,ruled,vlined]{algorithm2e}
\usepackage{booktabs}
\usepackage{diagbox}
\usepackage{color}

\expandafter\let\csname equation*\endcsname\relax
\expandafter\let\csname endequation*\endcsname\relax
\usepackage{amsmath}

\newtheorem{remark}{Remark}

\begin{document}
\title[ODE-DPS]{ODE-DPS: ODE-based Diffusion Posterior Sampling for Inverse Problems in Partial Differential Equation}
\author{Enze Jiang$^1$\footnotemark{}, Jishen Peng$^1$\footnotemark[\value{footnote}], \footnotetext{First Author and Second Author contribute equally to this work.}
Zheng Ma$^{1,2,3,4}$, Xiong-Bin Yan$^{1,4}$ \footnotemark{}}

\footnotetext{Corresponding author.}
\address{$^1$ School of Mathematical Sciences, Shanghai Jiao Tong University, Shanghai, China}
\address{$^2$ Institute of Natural Sciences, MOE-LSC, Shanghai Jiao Tong University, Shanghai, China}
\address{$^3$ Qing Yuan Research Institute, Shanghai Jiao Tong University, Shanghai, China}
\address{$^4$ CMA-Shanghai, Shanghai Jiao Tong University, Shanghai, China}
\ead{yanxb2015@163.com}

\begin{abstract}
  In recent years we have witnessed a growth in mathematics for deep learning, which has been used to solve inverse problems of partial differential equations (PDEs). However, most deep learning-based inversion methods either require paired data or necessitate retraining neural networks for modifications in the conditions of the inverse problem, significantly reducing the efficiency of inversion and limiting its applicability. To overcome this challenge, in this paper, leveraging the score-based generative diffusion model, we introduce a novel unsupervised inversion methodology tailored for solving inverse problems arising from PDEs. Our approach operates within the Bayesian inversion framework, treating the task of solving the posterior distribution as a conditional generation process achieved through solving a reverse-time stochastic differential equation. Furthermore, to enhance the accuracy of inversion results, we propose an ODE-based Diffusion Posterior Sampling inversion algorithm. The algorithm stems from the marginal probability density functions of two distinct forward generation processes that satisfy the same Fokker-Planck equation. Through a series of experiments involving various PDEs, we showcase the efficiency and robustness of our proposed method.
\end{abstract}

%
%
%
%
%

\section{Introduction}

Partial differential equations (PDEs) are fundamental concepts in various scientific and engineering fields, serving as the cornerstone of many disciplines. They are characterized by their dependence on parameters that govern their behavior, which allows them to describe a wide range of physical and operational scenarios. Solving partial differential equations under different parameter settings is meaningful because it enables predictions of physical laws under varying conditions. Various solution techniques, such as the finite element method (FEM), finite volumes, finite differences, and spectral methods, have been developed for tackling the problem \cite{quarteroni2008numerical,ern2004theory}. This process establishes a mapping from parameters to solutions, known as a \textit{parameter-to-solution map}.

Conversely, when given experimental data, there is often a need to infer the unknown parameters associated with the observed phenomenon, thereby determining the parameter regimes of the observed system. This situation is known as the \textit{inverse problem}. Inverse problems are generally more challenging than forward problems due to their ill-posed nature, where minor errors in observations can lead to significant inaccuracies in model parameters. This challenge is further exacerbated in scenarios involving observation noise, incomplete physics, and high-dimensional parameter spaces \cite{stuart2010inverse,vogel2002computational}.

So far, several regularization methods have been employed to address ill-posed inverse problems, including Tikhonov-type regularization methods, iterative regularization methods, and Bayesian inversion methods \cite{stuart2010inverse,vogel2002computational}. The Tikhonov regularization method stabilizes the solution of the inverse problem by incorporating regularization terms and carefully selecting regularization parameters \cite{kirsch2011introduction}. The iterative inversion method tackles the ill-posed nature of the inverse problem by implementing an iterative stopping criterion \cite{kirsch2011introduction}. The Bayesian inversion method transforms ill-posed inverse problems into well-posed posterior distribution problems through the introduction of suitable prior distributions \cite{stuart2010inverse}. Although significant progress has been witnessed in various inverse problems, they are hindered by their inability to effectively utilize data from past records. This limitation leads to inaccurate inversion, making it challenging for the resulting parameters to meet practical requirements.

Recently, deep learning has emerged as a potent tool for addressing inverse problems related to partial differential equations (PDEs). For example, Raissi et al. \cite{raissi2019physics} introduced physics-informed neural networks (PINNs) to tackle such problems. In this methodology, neural networks represent unknown parameters, leveraging physical equations to construct a loss function. Subsequently, stochastic gradient descent optimizes the neural network parameters, yielding numerical inversion results. In recent years, research in this field has flourished, incorporating various methodologies like weak adversarial networks for inverse problems \cite{bao2020numerical}, PDE-aware deep learning methods for inverse problems \cite{tenderini2022pde}, implicit neural representations for seismic inversion \cite{sun2023implicit}, and inverse obstacle scatter \cite{vlavsic2022implicit}. These approaches operate on a case-by-case basis, necessitating complete retraining for each new set of unknown parameters, measured data, or modifications to other conditions in the inverse problem.

To address the aforementioned case-by-case scenario, researchers have proposed inverse neural operators as a solution for PDE inverse problems. Fan et al. \cite{fan2020solving} introduces compact neural network architectures for regularized inverse operators to solve the EIT inverse problem. Molinaro et al. \cite{molinaro2023neural} propose a novel architecture called Neural Inverse Operators (NIOs) for solving PDE inverse problems. NIOs utilize a composition of DeepONets and FNOs to approximate mappings from operators to functions. Ovadia et al. \cite{ovadia2023vito} combine a U-Net-based architecture with a vision transformer to design an inverse neural operator for solving various PDE inverse problems. While these inverse neural operators have demonstrated effective numerical inversion results for PDE inverse problems, acquiring pairs of data remains a challenging issue.

Presently, generative models \cite{kawar2022denoising, song2021solving,baldi2012autoencoders,kingma2019introduction,creswell2018generative,whang2021composing} have emerged as an alternative approach for solving inverse problems, offering the advantage of not necessitating paired unknown parameters and measurement data. 
This method only requires partial prior data about the unknown parameters and does not rely on their corresponding measurement data. Generative models constitute a category of machine learning models designed to learn the underlying data distribution and produce new data samples resembling the original. Typical examples of generative models encompass Autoencoders \cite{baldi2012autoencoders}, Variational Autoencoders (VAEs) \cite{kingma2019introduction}, and Generative Adversarial Networks (GANs) \cite{creswell2018generative}, among others. In recent years, diffusion models have gained prominence as mainstream generative models. Simultaneously, they have demonstrated effectiveness in addressing inverse problems within the visual field, as evidenced by studies such as \cite{song2022pseudoinverse, moliner2023solving, kawar2022denoising,chung2023prompt}.

We present in this work
an \textit{ODE-based Diffusion Posterior Sampling (ODE-DPS)} algorithm for solving inverse problems in partial differential equations (PDEs). Specifically, by leveraging partial prior data of unknown parameters, combined with diffusion generative models and deep neural networks, we learn the prior distribution of the data, which is then applied in solving the Bayesian posterior distribution of the PDEs inverse problems. Then, we solve the posterior distribution using a conditional generation process
achieved by solving a reverse-time stochastic differential equation. To enhance the inversion accuracy, we propose the ODE-DPS regularization inversion algorithm. This algorithm originates from an equivalent form of a noising process for the score-based diffusion model, with its equivalence guaranteed by the marginal probability density functions satisfying the same Fokker-Planck equation as in the forward noising process. 
Furthermore, based on experimental observations, we propose modifying the original algorithm by replacing the $L_2$-norm data residual term with an adaptive norm residual term. The purpose is to reduce the errors near the boundaries in the inversion results obtained by the score-based diffusion sampling algorithm. We summarize the main novelties of this paper as follows:

\begin{itemize}
    \item We propose an ODE-based diffusion sampling algorithm to solve inverse problems of PDEs. This algorithm utilizes a limited amount of prior data and integrates score-based generative models to learn the prior distribution of unknown parameters. Subsequently, it is applied in a conditional distribution-based backward generation process induced by Bayesian inversion.
    \item By leveraging the Fokker-Planck equation satisfied by the marginal probability density function of the noising process in the generative model, we derived a new reverse ordinary differential equation. Subsequently, we discretized it numerically to develop our ODE-based DPS inversion algorithm.
    \item  Based on experimental observations, we propose modifying the diffusion posterior sampling algorithm by replacing the $L_2$-norm data residual term with an adaptive norm residual term to reduce the errors near the boundaries in the inversion results.
    \item ODE-DPS algorithm requires only a small amount of prior data for unknown parameters and does not require observation data, also referred to as labeled data, during the training phase. Once trained, this generative model can be applied to various inverse problem-solving tasks. Compared to traditional inversion algorithms, it enhances the inversion accuracy with minimal reduction in efficiency.
    \item Numerically, we compare our algorithm with traditional inversion methods across various inverse problems. It clearly shows that the ODE-DPS inversion algorithm can significantly improve the inversion accuracy.
\end{itemize}
The rest of this paper is organized as follows: In Section \ref{sec2}, we give an introduction to inverse problems. Then we elucidate the score-based diffusion model in Section \ref{sec3}. After that, we explain our method to solve the inverse problems in Section \ref{sec4}, apply it to some examples, and compare it with some traditional methods in Section \ref{sec5}. Finally, we end this paper with conclusions in Section \ref{sec6}.

\section{Preliminaries}\label{sec2}
In this section, we consider the inverse problem of a general form, aiming to estimate $\bar{x}_0\in D$ from the data represented by the operator equation:
\begin{equation}\label{2.1}
    \bar{y}_\delta=F(\bar{x}_0)+\xi_\delta,
\end{equation}
where $F:D\subset X\rightarrow Y$ is an operator between reflexive Banach spaces $(X,\|\cdot\|_X)$ and $(Y,\|\cdot\|_Y)$ with domain $D$. Here, $\bar{y}_{\delta}$ denotes the observed data with noise $\xi_{\delta}$.

 Generally speaking,  the inverse problem (\ref{2.1}) is ill-posed in the sense that the solution may not exist, be not unique, or not vary smoothly on the data due to the small noise caused by the measurement and recording process. Consequently, directly solving the operator equation (\ref{2.1}) for ill-posed inverse problems is not feasible. To tackle these challenges, regularization techniques are commonly employed. The Tikhonov regularization method \cite{tikhonov1977solutions, kirsch2011introduction} stands out as a well-established and effective approach for addressing inverse problems across various fields. This method replaces solving operator equation (\ref{2.1}) by minimizing the following problem:
     \begin{equation}\label{tik}
       \bar{x}_0:=\arg \min\limits_{\bar{x}}\| \bar{y}^\delta-F(\bar{x})\|_{Y}^2+\alpha\|\bar{x}\|_{X}^2,
    \end{equation}
where $\|\bar{x}\|_{X}$ denotes the regularization term, and $\alpha$ is a regularization parameter.
The efficacy of the Tikhonov regularization method in solving inverse problems is significantly influenced by the choice of regularization parameter and the regularization term norm $\|\cdot\|_{X}$. Optimal selection of a regularization parameter and regularization term can substantially enhance the accuracy of the inversion results. Nevertheless, determining the optimal parameter and norm are typically a challenging task. 

Additionally, another method widely employed in solving inverse problems is the iterative method \cite{ito2014inverse, bakushinsky2005iterative}. However, due to the ill-posed nature of the inverse problem, the iterative inversion method often encounters semi-convergence. This phenomenon manifests as a decrease followed by an increase in the relative error of the numerical inversion result with iteration process. Therefore, the success of this inversion method hinges on choosing a suitable stopping criterion. Morozov’s discrepancy principle \cite{nair2009morozov} is typically utilized for this purpose. Currently, integrating the discrepancy principle into iterative methods has emerged as another mainstream approach for addressing ill-posed inverse problems.

\begin{figure}[htbp]
    \centering
    \subfloat[The true source $f$.]{\label{fig_land_0}
        \includegraphics[width=0.3\columnwidth,height=.25\linewidth]{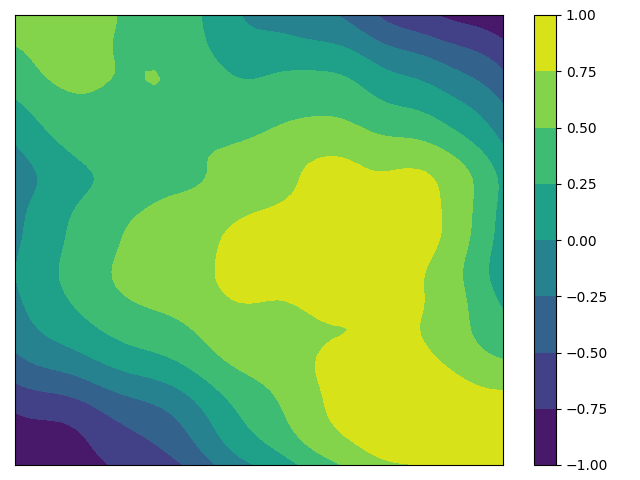}
    }
    \subfloat[The inversion $f$ by Landweber iteration method.]{\label{fig_land_a}
        \includegraphics[width=0.3\columnwidth,height=.25\linewidth]{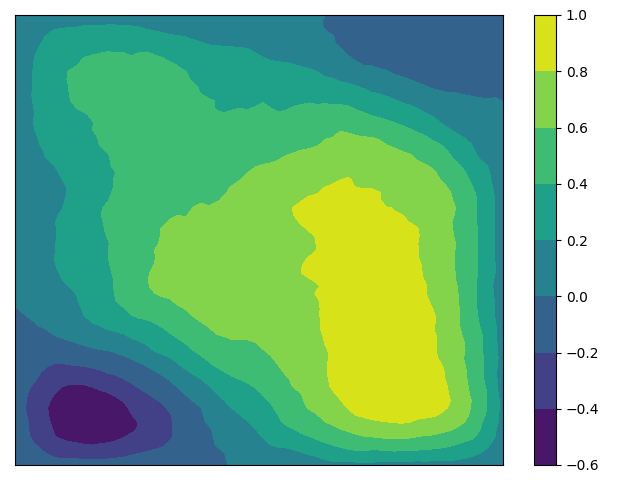}
    }
    \subfloat[The inversion $f$ by Tikhonov method.]{\label{fig_tik_b}
        \includegraphics[width=.3\columnwidth,height=.25\linewidth]{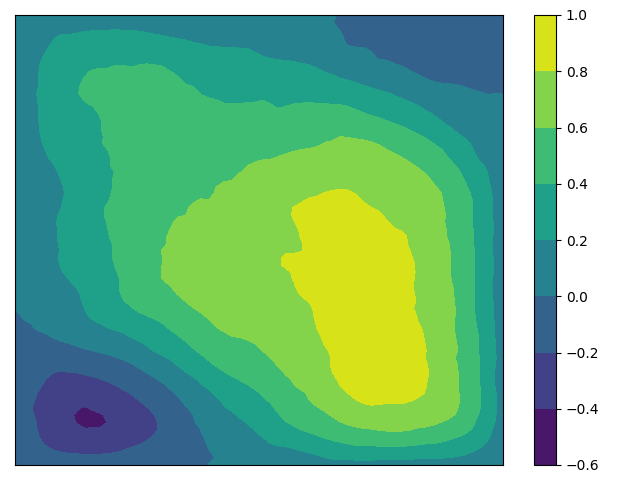}
    }

    \subfloat[The error by Landweber iteration method.]{\label{fig_land_c}
        \includegraphics[width=0.3\columnwidth,height=.25\linewidth]{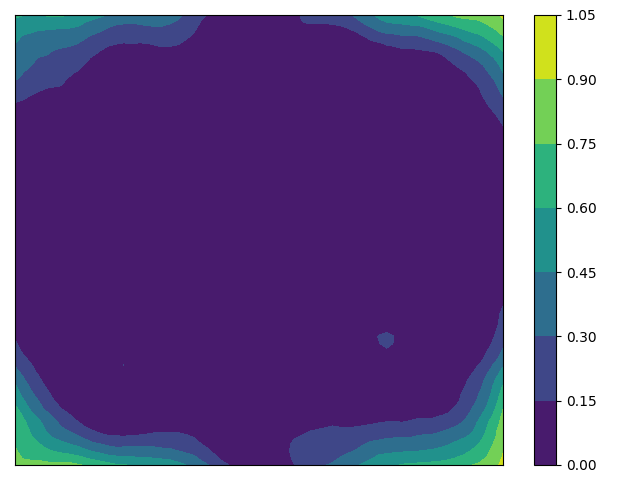}  
    } 
    \subfloat[The error by Tikhonov method.]{\label{fig_tik_d}
        \includegraphics[width=0.3\columnwidth,height=.25\linewidth]{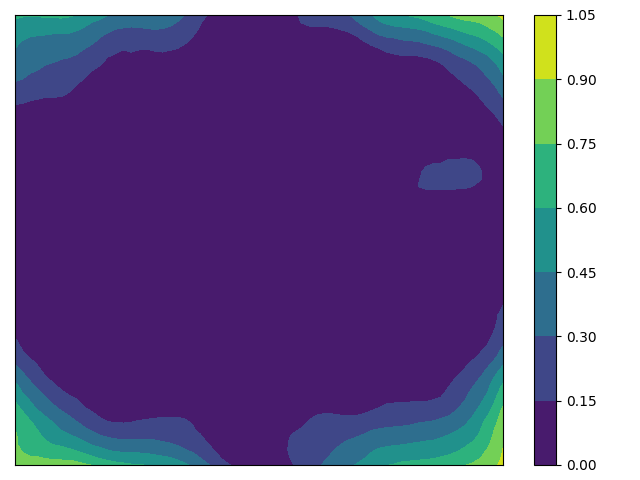}  
    }
    \caption{The inversion results of the problem (\ref{heq1}). (a) denotes the ground truth of source $f$. (b), (c) respectively provide the inversion results by Landweber iteration regularization and Tikhonov regularization method. (d), (e) respectively denotes the difference between true $f$ and the inversion source by Landweber iteration regularization and Tikhonov regularization method. The relative $l_2$ errors of the true and inversion source obtained by Landweber iteration and Tikhonov regularization are 34.4\% and 37.9\% respectively. Here, we set the regularization parameter $\alpha=0.005$ for Tikhonov regularization. }
    \label{fig_land}
\end{figure}

Although the inversion methods mentioned above have been applied to different inverse problems in the past, due to the ill-posedness of the inverse problem, we usually can only obtain low-precision inversion results. 
For example, let us consider the following equation:
\begin{equation}\label{heq1}
    \left\{\begin{aligned} 
        u_t(x,y,t) &= a\Delta u(x,y,t) + f(x,y),~(x,y,t)\in\Omega\times(0,T), \\  
        u(x,y,t) &= \psi(x,y,t),~(x,y,t)\in\partial\Omega\times (0,T),\\
        u(x,y,0)&= \phi(x,y), ~(x,y)\in\Omega,
    \end{aligned}\right.
\end{equation}
where $a$ is a positive constant. We assume that the initial value $\phi$ and boundary $\psi$ and the constant $a$ are known (Here we set $a=1, ~\psi(x,y,t)=0, ~\phi(x,y)=20\sin 2\pi x\sin 2\pi y,~\Omega=[0,1]\times[0,1], T=1$), and use the measurement data $u(x,y,T)$ at the terminal time $T=1$ to recover the unknown source term $f$. We set the time step size to 0.01 and the grid size to 1/63 to discretizate $\Omega$, and denote the measurement data with noise as $u^{\delta}(x,y, T)$. We apply Landweber iteration regularization \cite{hanke1995convergence} and Tikhonov regularization method \cite{kirsch2011introduction} for solving the ill-posed inverse source problem, and Morozov's discrepancy principle \cite{nair2009morozov} for stopping the Landweber iteration process. 
In the Tikhonov regularization method, we choose the norm of the regularization term to be the $l_2$ norm, with a regularization parameter of $\alpha=0.005$.
The inversion results are shown in Figure \ref{fig_land}. Compared with the exact solution in Figure \ref{fig_land_0}, the numerical results of the inversion in Figures \ref{fig_land_a}, \ref{fig_tik_b} only captures large-scale features, and small-scale features are difficult to capture. In addition, the errors in Figures \ref{fig_land_c}, \ref{fig_tik_d} also reflect that Landweber iteration regularization method and Tikhonov regularization method are difficult to give accurate inversion results. 
The main reasons why this inversion method struggles to achieve high-precision inversion results are primarily twofold:
\begin{enumerate}
    \item Although we only consider a linear inverse problem in this numerical example, their inherent ill-posed nature persists. Hence, it remains challenging to obtain an accurate solution for the ill-posed problem. 
    \item While we have utilized the mainstream regularization methods to solve the aforementioned inverse problem, we have not incorporated any prior information about the inverse problem into the inversion process, apart from selecting the initial value and the norm of regularization term. In reality, leveraging prior information past-recorded of inverse problems can significantly enhance inversion accuracy.
\end{enumerate}

Therefore, it is necessary for us to design a new regularized inversion algorithm to solve the ill-posed inverse problem. 
In recent years, deep learning-based inversion methods, as another form of regularization, have been applied across various fields of inverse problems, showing a trend of outperforming traditional inversion methods. Their advantage primarily stems from their ability to learn prior information about inverse problem solutions from previously recorded data, which can then be applied to the inversion process. Building upon this notion, this paper proposes a novel regularization inversion method. This method leverages partial prior information of unknown parameters, deep neural networks, score-based diffusion generative model, and the traditional Bayesian inversion framework to introduce a completely new regularization inversion approach. The following section will present the interpretation of our inversion method.

\section{Score-Based Diffusion Model}\label{sec3}
In this article, our main focus is on constructing new regularization algorithms using generative models and Bayesian inversion methods. Deep generative models are widely used in many subfields of AI and Machine Learning. These models, when combined with stochastic optimization methods, excel at capturing complex, high-dimensional data distributions such as images, text, and speech. Today, popular deep generative models include variational autoencoders (VAEs), generative adversarial networks (GANs), auto-regressive models, and so on. In particular, diffusion models, also known as diffusion probabilistic models or score-based generative models, represent a class of latent variable generative models that are found in applications for diverse inverse problem tasks, including image denoising \cite{xie2023diffusion}, inpainting \cite{corneanu2024latentpaint}, super-resolution \cite{yue2024resshift}, and medical imaging \cite{song2021solving}. These models have demonstrated state-of-the-art performance in these fields. Consequently, our objective in this study is to devise a novel regularization inversion method for addressing PDE inverse problems by integrating score-based diffusion models with traditional Bayesian inversion methods.

Score-based generative models represent a continuum of distributions evolving over time through a diffusion process. This process gradually transforms a data point into random noise, as described by a predetermined stochastic differential equation (SDE) independent of the data and devoid of trainable parameters. Reversing this process enables us to smoothly convert random noise into data for sample generation. Importantly, this reverse process adheres to a reverse-time SDE, derivable from the forward SDE based on the score of marginal probability densities over time. Consequently, we can approximate the reverse-time SDE by training a time-dependent neural network to estimate the scores, facilitating sample production using numerical SDE solvers. Below we will explain the score-based generative model in detail.

\subsection{Perturbing Data with SDE}
Score-based diffusion model \cite{song2020score} defines the generative process as the \textit{reverse} of the noising process. Specifically, we construct a diffusion process $\mathbf{x}(t)$ indexed by a continuous time variable $t\in [0,\hat{T}]$, such that $\mathbf{x(0)}\sim p(x(0))=p_0$, for which
    we have a datasets of i.i.d samples, and $\mathbf{x}(\hat{T})\sim p(x(\hat{T}))=p_{\hat{T}}$, for which we have a tractable form to generative samples efficiently. In other words, $p_0$ is the data distribution and $p_{\hat{T}}$ is the prior distribution. The diffusion process
can be modeled as the solution to an $\mathrm{It\hat{o}}$  stochastic different equation (SDE):
\begin{equation}\label{2.2}
    d\mathbf{x}=-\frac{\beta(t)}{2}\mathbf{x}dt+\sqrt{\beta(t)}d\mathbf{\omega}
\end{equation}
where $\beta(t)>0$ is the noise schedule of the process, typically taken to be monotonically increasing linear function of  $t$ \cite{ho2020denoising}, and 
$\mathbf{\omega}$ is the standard Wiener process (a.k.a. Brownian motion). The SDE has a 
unique strong solution as long as the coefficients are globally Lipschitz in
both state and time \cite{oksendal2003stochastic}. We hereafter denote
by $p(\mathbf{x}(t))$ the probability density of $\mathbf{x}(t)$, and use $p(\mathbf{x}(t)|\mathbf{x}(s))$
to denote the transition kernel from $\mathbf{x}(s)$ to $\mathbf{x}(t)$,  where $0\le s<t\le \hat{T} $.  

\subsection{Generative Samples by Reversing the SDE}
The generative process is to recover the data distribution, i.e, obtain samples $\mathbf{x}(0)\sim p(x(0))$,  by starting from samples of $\mathbf{x}(\hat{T})\sim p(x(\hat{T}))$ and reversing the diffusion process
(\ref{2.2}). A remarkable result from Anderson \cite{anderson1982reverse} states that the \textit{reverse} of a diffusion process is also a diffusion process, running backwards in time and given by the reverse-time SDE:
\begin{equation}
    d\mathbf{x}=[-\frac{\beta(t)}{2}\mathbf{x}-\beta(t)\nabla_{\mathbf{x}(t)}\log p(\mathbf{x}(t))]dt+\sqrt{\beta(t)}d\bar{\mathbf{\omega}},
    \label{2.3}
\end{equation}
where $\bar{\mathbf{\omega}}$ is a standard Wiener process when time flows backwards from $\hat{T}$ to $0$, and $dt$ is an infinitesimal negative timestep. Once the score of each
marginal distribution, $\nabla_{\mathbf{x}(t)}\log p(\mathbf{x}(t))$, is known for all $t$,  we can derive the reverse diffusion process from equation (\ref{2.3}) and simulate it to obtain
samples from distribution $p_0$.

\subsection{Estimating Score for the SDE}
The score $\nabla_{\mathbf{x}(t)}\log p(\mathbf{x}(t))$ can be estimated by training a score-based model on samples, i.e, we can train a time-dependent score-based model $s_\theta(\mathbf{x},t)$ to approximate the score function $\nabla_{\mathbf{x}(t)}\log p(\mathbf{x}(t))$ with denoising score matching (DSM) \cite{vincent2011connection}. The optimization objective was proved equivalent to the following:
\begin{equation}
    \theta^*=\mathop{\arg\min}\limits_\theta\mathbb{E}_{t\sim U(0, \hat{T}),~\mathbf{x}(t)\sim p(\mathbf{x}(t)|\mathbf{x}(0)),~\mathbf{x}(0)\sim p(\mathbf{x}(0))}[\|s_\theta(\mathbf{x}(t),t)-\nabla_{\mathbf{x}(t)}\log p(\mathbf{x}(t)|\mathbf{x}(0))\|_2^2],
    \label{eq5}
\end{equation}
where $\epsilon\simeq 0$  is a small positive constant. With sufficient data and model capacity, score matching ensures that the optimal solution to equation (\ref{eq5}) denoted by $s_{\theta^{*}}(\mathbf{x}(t), t)$.  Then, one can use the approximation 
$\nabla_{\mathbf{x}(t)}\log p(\mathbf{x}(t))\simeq s_{\theta^{*}}(\mathbf{x}(t), t)$ as a plug in estimate to replace the score function in (\ref{2.3}).  Discretization of (\ref{2.3}) and solving using, e.g. Euler-Maruyama discretization, amounts to sampling from the data distribution $p_0$, the goal of generative modeling.

\section{Solving inverse problem by Score-based generative model}\label{sec4}
Consider the inverse problem mentioned in (\ref{2.1}), and  we assume that  the unknown variable $\bar{x}_0$ and the observed data $\bar{y}_{\delta}$ are random variables, $\xi_{\delta}$ is the measurement noise. In the case of white Gaussian noise, $\xi_{\delta}\sim \mathcal{N}(0,\sigma^2 I)$. In the Bayesian framework, one utilizes $p(\bar{x}_0)$ as the prior,
and denote the likelihood $p(\bar{y}_{\delta}|\bar{x}_0)\sim \mathcal{N}(\bar{y}_{\delta}|F(\bar{x}_0), \sigma^2 I)$.
By Bayesian rules', the posterior distribution $p(\bar{x}_0|\bar{y}_{\delta})$ 
denoted by 
\begin{equation*}\label{pd_01}
p(\bar{x}_0|\bar{y}_{\delta})=\frac{1}{p(\bar{y}_{\delta})}p(\bar{y}_{\delta}|\bar{x}_0)p(\bar{x}_0).
\end{equation*}

Then the inverse problem (\ref{2.1}) is transformed into the problem of solving the posterior distribution $p(\bar{x}_0|\bar{y}_{\delta})$. To do it, we use the forward diffusion process (\ref{2.2}) to obtain a family of diffused distribution $p(\bar{x}(t)|\bar{y}_{\delta})$ with the initial distribution $p(\bar{x}_0|\bar{y}_{\delta})$ and apply Anderson's Theorem to derive the conditioned reverse time SDE
\begin{equation}\label{eq7}
    d\bar{x}(t)=[-\frac{\beta(t)}{2}\bar{x}(t)-\beta(t)\nabla_{\bar{x}(t)}
    \log {p(\bar{x}(t)|\bar{y}_{\delta})}]dt+\sqrt{\beta(t)}d\bar{\mathbf{\omega}}.
\end{equation}
Then the solution to the inverse problem, that is, solving the posterior distribution $p(\bar{x}_0|\bar{y}_{\delta})$, is transformed into sampling from the conditional reverse time SDE (\ref{eq7}).
Since $p(\bar{x}(t)|\bar{y}_{\delta})\propto p(\bar{x}(t))p(\bar{y}_{\delta}|\bar{x}(t))$, the score $\nabla_{\bar{x}(t)}\log p(\bar{x}(t)|\bar{y}_{\delta})$ in (\ref{eq7})
can be computed easily by 
\begin{equation}\label{2.4}
\nabla_{\bar{x}(t)}\log p(\bar{x}(t)|\bar{y}_{\delta})=\nabla_{\bar{x}(t)}\log p(\bar{x}(t))+\nabla_{\bar{x}(t)}\log p(\bar{y}_{\delta}|\bar{x}(t)).
\end{equation}
Combing (\ref{eq7}) with (\ref{2.4}), we have 
\begin{equation}
    d\bar{x}(t)=[-\frac{\beta(t)}{2}\bar{x}(t)-\beta(t)(\nabla_{\bar{x}(t)}\log p(\bar{x}(t))+\nabla_{\bar{x}(t)}\log p(\bar{y}_{\delta}|\bar{x}(t)))]dt+\sqrt{\beta(t)}d\bar{\mathbf{\omega}}.
    \label{2.5}
\end{equation}
Leveraging (\ref{2.5}), scored-based generative models provide a way 
to solve the Bayesian inverse problem, i.e., to solve the posterior distribution
$p(\bar{x}_0|\bar{y}_{\delta})$.  In (\ref{2.5}), we have two terms that should be computed: the score
function $\nabla_{\bar{x}(t)}\log p(\bar{x}(t))$ and the likelihood function $\nabla_{\bar{x}(t)}\log p(y_{\delta}|\bar{x}(t))$.
To compute the former term $\nabla_{\bar{x}(t)}\log p(\bar{x}(t))$, we can simply use the pre-trained score function $s_{\theta^{*}}(\bar{x}(t),t)$. However, the latter term is hard to acquire 
in closed-form due to the dependence on the time $t$, as there only exists dependence between $\bar{y}_{\delta}$ and $\bar{x}_0$.

\subsection{Approximate the Likelihood}
In this part,  we try to approximate $\log p(\bar{y}_{\delta}|\bar{x}(t))$. Since $\bar{x}(t)$ and $\bar{y}_{\delta}$ are conditionally independent on $\bar{x}_0$, we have 
\begin{equation}
\begin{aligned}
    p(\bar{y}_{\delta}|\bar{x}(t))&=\int p(\bar{y}_{\delta}|\bar{x}_0,\bar{x}(t))p(\bar{x}_0|\bar{x}(t))d\bar{x}_0\\
    &=\int p(\bar{y}_{\delta}|\bar{x}_0)p(\bar{x}_0|\bar{x}(t))d\bar{x}_0\\
    &=\mathbb{E}_{\bar{x}_0\sim p(\bar{x}_0|\bar{x}(t))}[p(\bar{y}_{\delta}|\bar{x}_0)].
\end{aligned}
\end{equation}
Chung et al \cite{chung2022diffusion} suggests the following approximation:
\begin{equation*}
    p(\bar{y}_{\delta}|\bar{x}(t))\simeq p(\bar{y}_{\delta}|\hat{\bar{x}}_0),
\end{equation*}
where 
\begin{equation*}
\hat{\bar{x}}_0:=\mathbb{E}[\bar{x}_0|\bar{x}(t)]=\frac{1}{\sqrt{\bar{\alpha}(t)}}(\bar{x}(t)+(1-\bar{\alpha}(t))\nabla_{\bar{x}(t)}\log p(\bar{x}(t))).
\end{equation*}
Replacing $\nabla_{\bar{x}(t)}\log p(\bar{x}(t))$ with the trained neural network
$s_{\theta^*}(\bar{x}(t),t))$, we have the following approximation:
\begin{equation*}
    \hat{\bar{x}}_0\simeq\frac{1}{\sqrt{\bar\alpha(t)}}(\bar{x}(t)+(1-\bar\alpha(t))s_{\theta^*}(\bar{x}(t),t)),
\end{equation*}

Now,  we turn to the calculation of $p(\bar{y}_{\delta}|\hat{\bar{x}}_0)$. The likely function takes the form 
\begin{equation*}
    p(\bar{y}_{\delta}|\hat{\bar{x}}_0)=\frac{1}{\sqrt{(2\pi)^n\sigma^{2n}}}\exp[-\frac{\|\bar{y}_{\delta}-F(\hat{\bar{x}}_0)\|_2^2}{2\sigma^2}],
\end{equation*}
where $n$ is the dimension of $\bar{y}_{\delta}$. We can calculate an approximation value of $\nabla_{\bar{x}(t)}\log p(\bar{y}_{\delta}|\bar{x}(t))$ in (\ref{2.5}) by 
\begin{equation}
    \nabla_{\bar{x}(t)}\log p(\bar{y}_{\delta}|\bar{x}(t))\simeq
    \nabla_{\mathbf{x}(t)}\log p(\bar{y}_{\delta}|\hat{\bar{x}}_0)\simeq -\frac{1}{\sigma^2}\nabla_{\bar{x}(t)}\|\bar{y}_{\delta}-F(\hat{\bar{x}}_0(\bar{x}(t)))\|_2^2.
    \label{eq15}
\end{equation}
Combining the equation (\ref{eq15}) and (\ref{2.4}), we can get the approximation  
\begin{equation}\label{eq16}
     \nabla_{\bar{x}(t)}\log p(\bar{x}(t)|\bar{y}_{\delta})\simeq s_{\theta*}(\bar{x}(t),t)-\zeta\nabla_{\bar{x}(t)}\|\bar{y}_{\delta}-F(\hat{\bar{x}}_0(\bar{x}(t)))\|_2^2,
\end{equation}
$\zeta:=\frac{1}{\sigma^2}$ is the step size, where we replace
$\nabla_{\bar{x}(t)}\log p(\bar{x}(t))$ with the trained neural network $s_{\theta^*}(\bar{x}(t),t))$. 

Through the above discussion, and  incorporateing (\ref{eq16}) into (\ref{2.5}) and using the usual ancestral sampling, we can get the Diffusion Posterior Sampling (DPS) algorithm \ref{al11}.  

\begin{algorithm}[t]
    \textbf{Require}: $N,\bar{y}_{\delta},s_{\theta^{*}},\zeta, \{\beta_i\}_{i=1}^{N}$; \\
    \textbf{Result}: $\bar{x}_0$; \\
    $\bar{x}_{N}\sim\mathcal{N}(0,\mathbf{I})$; \\
    $\alpha_i := 1-\beta_i,\bar{\alpha_i}:=\prod_{j=1}^{i}\alpha_i$\\
    \For{$i=N$ \KwTo $0$}{
        $\hat{s} = s_{\theta^{*}}(\bar{x}_i, i)$ \\
        $\hat{\bar{x}}_0=\frac{1}{\sqrt{\bar{\alpha}_i}}(\bar{x}_i+(1-\bar{\alpha}_i)\hat{s})$\\
        $z\sim\mathcal{N}(0,I)$\\
        $\bar{x}_{i-1}'=\frac{\sqrt{\alpha_i}(1-\bar{\alpha}_{i-1})}{1-\bar{\alpha}_i}\bar{x}_i+\frac{\sqrt{\bar{\alpha}_{i-1}}\beta_i}{1-\bar{\alpha}_i}\hat{\bar{x}}_0+\tilde{\sigma}_i z$\\
        $\bar{x}_{i-1} = \bar{x}_{i-1}'-\zeta \nabla_{\bar{x}_i}\|\bar{y}_{\delta}-F(\hat{\bar{x}}_0)\|_{2}^2$. \\
    }
    \caption{Diffusion Posterior Sampling (DPS) for Inverse problem.}
    \label{al11}
\end{algorithm}
\begin{figure}[t]
    \centering
    \subfloat[The true $f$.]{\label{fig_DPS_a}
        \includegraphics[width=0.32\columnwidth,height=.28\linewidth]{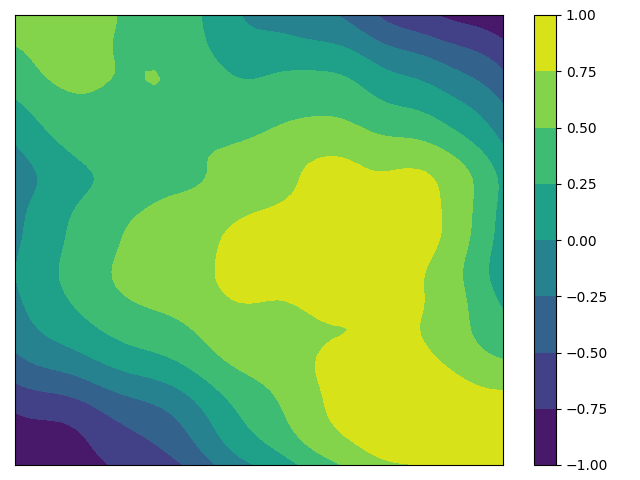}
    }
    \subfloat[The inversion $f$.]{\label{fig_DPS_b}
        \includegraphics[width=.32\columnwidth,height=.28\linewidth]{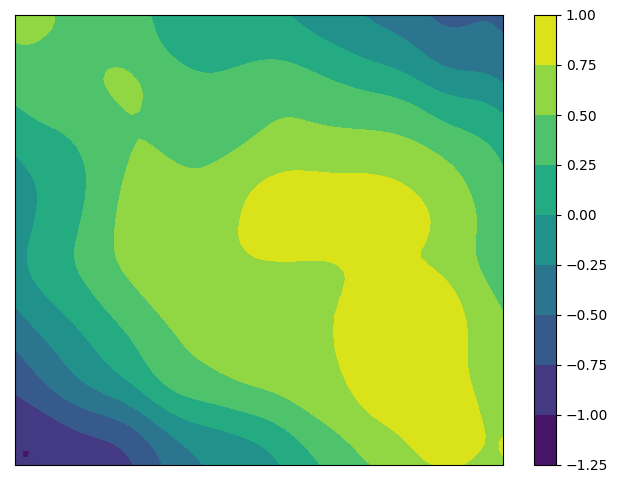}
    }
    \subfloat[Difference between true source and inversion source.]{\label{fig_DPS_c}
        \includegraphics[width=0.32\columnwidth,height=.28\linewidth]{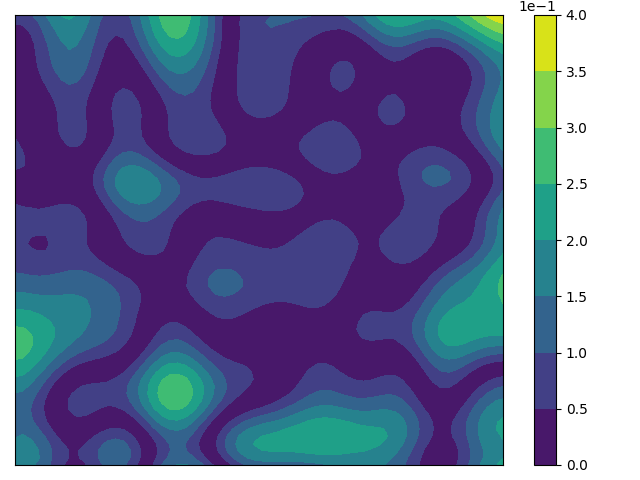}
    }
    \caption{The inversion results of the problem (\ref{heq1}). (a) denotes the ground truth of source $f$. (b) denotes the inversion results by Diffusion Posterior Sampling algorithm \ref{al11}. (c) denotes the difference between true $f$ and the inversion source. The relative $l_2$ error of the true and inversion source is  18.4\%. In algorithm \ref{al11}, we set $N=1000$.}
    \label{fig_DPS}
\end{figure}

In Figure \ref{fig_DPS}, we present the numerical results obtained using the DPS inversion algorithm to solve the ill-posed inverse problem described in Section \ref{sec2}. The network architecture, parameters, and training details of the neural network $s_{\theta^{*}}(\bar{x}(t),t)$ will be deferred and presented later in the paper.
It is evident from Figure \ref{fig_DPS_b} that the DPS inversion algorithm significantly improves the inversion performance compared to the Landweber iteration and Tikhonov regularization method. The relative $l_2$ error of the inversion results obtained by the DPS algorithm is 18.4\%, as opposed to the relative $l_2$ error 34.4\% and 37.9\% produced by the Landweber iteration and Tikhonov regularization respectively, indicating a substantial enhancement in the inversion accuracy with the DPS algorithm. However, in Figure \ref{fig_DPS_c}, we observe larger errors near the boundaries of the region. Therefore, reducing the errors near the boundaries of the inversion results can improve the solution accuracy. 
Moreover, in the experiment, we find that the original DPS algorithm is unstable, which means the inversion results we get exist some fluctuations due to the randomness of the algorithm. Therefore, to achieve higher accuracy and better stability in the inversion algorithm, we need to modify the original DPS inversion algorithm.
Based on this, in the following, we will provide an improved version of the DPS algorithm.

\subsection{ODE-based Diffusion Posterior Sampling for inverse problems}

We reconsider a stochastic differential equation in the following form:
\begin{equation*}
    d\mathbf{x}(t)=f(\mathbf{x}(t),t)dt+g(t)dw,
\end{equation*}
which is more generative than SDE (\ref{2.2}). Suppose that $p(\mathbf{x}(t))$ is the marginal probability density function, then the corresponding Fokker-Planck equation \cite{song2020score} is:
\begin{equation}    \label{FP1}
    \frac{\partial p(\mathbf{x}(t))}{\partial t}=-\nabla_{\mathbf{x}(t)}\cdot[f(\mathbf{x}(t),t)p(\mathbf{x}(t))]+\frac{1}{2}g^2(t)\nabla_{\mathbf{x}(t)}\nabla_{\mathbf{x}(t)}p(\mathbf{x}(t)).
\end{equation}
Using the fact that $\nabla_{\mathbf{x}(t)}(\log p(\mathbf{x}(t)))p(\mathbf{x}(t))=\nabla_{\mathbf{x}(t)}p(\mathbf{x}(t))$, the equation (\ref{FP1}) can be rewritten into the following form:
\begin{equation}
\begin{aligned}
    \frac{\partial p(\mathbf{x}(t))}{\partial t}=&-\nabla_{\mathbf{x}(t)}\cdot[f(\mathbf{x}(t),t)p(\mathbf{x}(t))-\frac{1}{2}(g^2(t)-\mu^2(t))\nabla_{\mathbf{x}(t)}p(\mathbf{x}(t))] \\
    &+\frac{1}{2}\mu^2(t)\nabla_{\mathbf{x}(t)}\nabla_{\mathbf{x}(t)}p(\mathbf{x}(t))\\
    =&-\nabla_{\mathbf{x}(t)}\cdot[(f(\mathbf{x}(t),t)-\frac{1}{2}(g^2(t)-\mu^2(t))\nabla_{\mathbf{x}(t)}\log p(\mathbf{x}(t)))p(\mathbf{x}(t))] \\
    &+\frac{1}{2}\mu^2(t)\nabla_{\mathbf{x}(t)}\nabla_{\mathbf{x}(t)}p(\mathbf{x}(t)),
    \label{FP2}
\end{aligned}
\end{equation}
where $\mu(t)$ is an arbitrary smooth function. By comparing (\ref{FP1}) and (\ref{FP2}), the SDE corresponding to equation (\ref{FP2}) is 
\begin{equation*}
   d\mathbf{x}(t)=[f(\mathbf{x}(t),t)-\frac{1}{2}(g^2(t)-\mu^2(t))\nabla_{\mathbf{x}(t)}\log p(\mathbf{x}(t))]dt+\mu(t)dw.
\end{equation*}
Again, by Anderson \cite{anderson1982reverse}, we obtain the reverse SDE of the above SDE as follows:
\begin{equation*}
    d\mathbf{x}(t)=[f(\mathbf{x}(t),t)-\frac{1}{2}(g^2(t)+\mu^2(t))\nabla_{\mathbf{x}(t)}\log p(\mathbf{x}(t))]dt+\mu(t)dw.
\end{equation*}
Let $\mu(t)=0, f(\mathbf{x}(t),t)=-\frac{\beta(t)}{2}\mathbf{x}, ~g(t)=\sqrt{\beta(t)}$, we can get the following ordinary differential equation,
\begin{equation}\label{4.2}
    d\mathbf{x}(t) = [-\frac{\beta(t)}{2}\mathbf{x}(t)-\frac{\beta(t)}{2}\nabla_{\mathbf{x}(t)}\log p(\mathbf{x}(t))]dt.
\end{equation}

By taking a negative time step, the differential equation (\ref{4.2}) defines a new generative process from $x(\hat{T})\sim p(x(\hat{T}))$. 
Similar to the previous discussion, for solving the inverse problem (\ref{2.1}), we only need to solve the following reverse ordinary differential equation:
\begin{equation}\label{eqODE}
    d\bar{x}(t) = [-\frac{\beta(t)}{2}\bar{x}(t)-\frac{\beta(t)}{2}(\nabla_{\bar{x}(t)}\log p(\bar{x}(t))+\nabla_{\bar{x}(t)}\log p(\bar{y}_{\delta}|\bar{x}(t)))]dt.
\end{equation}

Additionally, we observed in experiments that the gradient values near the boundaries in the gradient descent steps of the DPS inversion algorithm are relatively small. This suggests that the values near the boundaries remain nearly unchanged when updating the parameters for inversion, leading to significant errors near the boundaries in the inversion results. To mitigate this issue, we propose modifying the $L_2$ norm of the data residual term to an adaptive norm. Specifically, we introduce a new weighted norm for the function (or matrix) $M$, defined as follows:
 \begin{equation*}
\begin{aligned}
    \|M(x)\|_\rho &= \int_\Omega \rho(x)M(x)^2 dx,\\
    \|M\|_\rho &= \sum_i\sum_j \rho_{ij}M_{ij}^2,
\end{aligned}
\end{equation*}
where $\Omega$ is the domain of the inversion parameters. $\rho(x)$($\{\rho_{ij}\}$) is a weighting function(weighting matrix), which is set as:
\begin{equation}\label{de_rho}
\begin{aligned}
    \rho(x) &= \frac{(\max_{x \in \Omega}\frac{\partial \|F(\bar{x})\|_2}{\partial \bar{x}}|_{\bar{x}=\hat{\bar{x}}_0}(x))^2+c}{(\frac{\partial \|F(\bar{x})\|_2}{\partial \bar{x}}|_{\bar{x}=\hat{\bar{x}}_0}(x))^2+c}\\
    \rho_{ij} &= \frac{(\max_{x_{ij} \in \Omega}\frac{\partial \|F(\bar{x})\|}{\partial \bar{x}_{ij}}|_{\bar{x}=\hat{\bar{x}}_0})^2+c}{(\frac{\partial \|F(\bar{x})\|_2}{\partial_{\bar{x}_{ij}}} |_{\bar{x}=\hat{\bar{x}}_0})^2+c},
\end{aligned}
\end{equation}
where $c$ is a small constant to guarantee the stability of gradients. Therefore, we modify equation (\ref{eq16}) to the following form
\begin{equation}
     \nabla_{\bar{x}(t)}\log p(\bar{x}(t)|\bar{y}_{\delta})\simeq s_{\theta*}(\bar{x}(t),t)-\zeta_i\nabla_{\bar{x}(t)}\|\bar{y}_{\delta}-F(\hat{\bar{x}}_0(\bar{x}(t)))\|_\rho^2.
    \label{eq19}
\end{equation}
Based on the above discussion, and combing formulas (\ref{eqODE}) with (\ref{eq19}), and using the ancestral sampling, we get an ODE-based DPS (ODE-DPS) algorithm as summarized in algorithm \ref{al2}. 
\begin{algorithm}[t]
    \textbf{Require}: $N,~\bar{y}_{\delta},~s_{\theta^{*}},~\zeta,~\gamma,~c, \{\beta_i\}_{i=1}^N$;\\
    \textbf{Result}: $\bar{x}_0$;\\
    $\alpha_i := 1-\beta_i,\bar{\alpha_i}:=\prod_{j=1}^{i}\alpha_i$\\
    $\bar{x}_N\sim\mathcal{N}(0,\mathbf{I})$; \\
    \For{$i=N$ \KwTo $0$}{
        $\hat{s} = s_{\theta^{*}}(\bar{x}_i, i)$ \\
        $\hat{\bar{x}}_0=\frac{1}{\sqrt{\bar{\alpha}_i}}(\bar{x}_i+\frac{1}{2}(1-\bar{\alpha}_i)\hat{s})$\\
$\bar{x}_{i-1}'=\frac{\sqrt{\alpha_i}(1-\bar{\alpha}_{i-1})}{1-\bar{\alpha}_i}\bar{x}_i+\frac{\sqrt{\bar{\alpha}_{i-1}}\beta_i}{1-\bar{\alpha}_i}\hat{\bar{x}}_0$ \\
        $\zeta_i=\zeta\gamma^{[\frac{N-1}{100}]}$ \\
        $\rho(x) = \frac{(\max_{x \in \Omega}\frac{\partial \|F(\bar{x})\|_2}{\partial \bar{x}}|_{\bar{x}=\hat{\bar{x}}_0}(x))^2+c}{(\frac{\partial \|F(\bar{x})\|_2}{\partial \bar{x}}|_{\bar{x}=\hat{\bar{x}}_0}(x))^2+c}$ \\
        $\bar{x}_{i-1} = \bar{x}_{i-1}'-\zeta_i \nabla_{\bar{x}_i}\|\bar{y}_{\delta}-F(\hat{\bar{x}}_0)\|_{\rho}^2$. \\
    }
    \caption{ODE-DPS for Inverse problem}
    \label{al2}
\end{algorithm}

\begin{remark}
It is worth noting that, in Algorithm \ref{al2}, there is no need to retrain the neural network $s_{\theta}(\bar{x},t)$. This approach requires training the neural network only once on the dataset, making it applicable to various inverse problem-solving tasks. This significantly reduces the inversion time compared to the inverse methods based on Physics-Informed Neural Networks (PINNs) \cite{raissi2019physics}, thus enhancing the inversion efficiency. Moreover, our inversion algorithm, akin to traditional regularization methods, incorporates the model information of the inverse problem. Compared to neural operator-based inversion methods \cite{fan2020solving, ovadia2023vito}, it improves generalization and enhances the algorithm's interpretability.
\end{remark}

\subsection{U-Net Architecture}
\begin{figure}
    \centering
    \includegraphics[width=\columnwidth,height=0.5\linewidth]{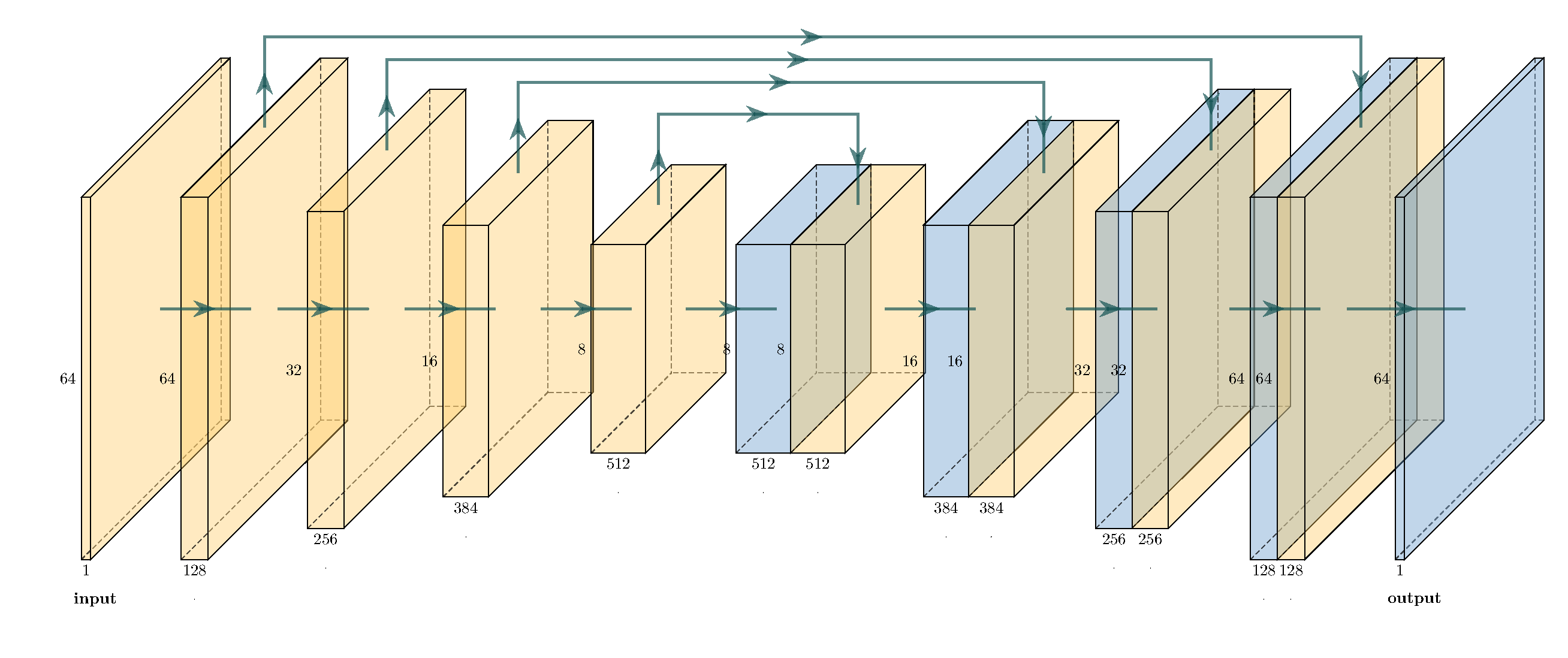}
    \caption{The Simplified U-Net model structure.}
    \label{fig:Unet}
\end{figure}

For diffusion models, U-Net architecture \cite{ronneberger2015u} is always used. The Simplified U-Net model structure is shown in Figure \ref{fig:Unet}, and we specifically illustrate the architecture of the U-Net model along with the parameters we utilized in \ref{B}.
\par The U-Net model consists of a downsampling path and an upsampling path. These two path are connected with skip connections that do not change the spatial size. Downsampling path employs numerous convolutional layers, decreasing the spatial size and increase the number of feature channels. Upsampling path increases the spatial size and decreasing the number of feature channels. For the skip connections, padding in convolutional layers of the downsampling path ensures that the tensors we concatenate have the same spatial size. Multi-heads attention layers \cite{vaswani2017attention} are also used in this work, which is found that can improve the behavior of the U-Net model. In addtion, we need to add a projection of the timestep embedding into each residual block. 

\section{Experiment}\label{sec5}
In this section, we show the performance of our method by three different inverse problems. Experimental performance is evaluated using both qualitative analysis and
quantitative metric. We define the relative $l_2$ error as follows:
\begin{equation*}
    \textrm{relative}~l_2~\textrm{error}=\frac{\|\bar{x}_{ref}-\bar{x}_{rec}\|}{\|\bar{x}_{ref}\|},
\end{equation*}
where $\bar{x}_{ref}$ is the ground truth and $\bar{x}_{rec}$ is the reconstructed parameter.
The relative $l_2$ error measures the degree of accuracy in reconstructing the unknown parameters. 

To show the advantages of our proposed inversion algorithm, we will contrast it with the traditional regularization inversion methods, namely the Landweber iterative regularization method \cite{hanke1995convergence, kirsch2011introduction} and the Tikhonov regularization method \cite{kirsch2011introduction}, in the paper. We illustrate the methods as follows:
\begin{enumerate}
    \label{traditional}
    \item \textbf{Method 1}: Landweber iteration regularization. It is a widely used regularization inversion method applied in solving ill-posed inverse problems, and various extensions of this method have been developed for different inverse problems. In this paper, we only focus on the most basic version \cite{hanke1995convergence}, i.e., given $\bar{x}_0$ is an initial guess which may incorporate a prior knowledge of an exact solution for inverse problem, the iteration is defined via the adjoint $F^{'}(\cdot)^{*}$ of $F^{'}(\cdot)$ as follows:
    \begin{equation*}
        \bar{x}_{k+1}^{\delta}=\bar{x}_{k}^{\delta}+\zeta_k F^{'}(\bar{x}_k^{\delta})^{*}(\bar{y}^{\delta}-F(\bar{x}_k^{\delta})),~k=0,1,2,\cdots,
    \end{equation*}
    where $F^{'}(\cdot)$ is the Fr\'echet derivative of forward operator $F(\cdot)$ and $\zeta_k$ is the step size of $k$th step. To obtain a stable approximation of the method, the iteration 
    process need to be stopped after $k_{*}=k_{*}(\delta)$ steps according to a 
    generalized discrepancy principle, i.e.,
    \begin{equation*}
    \|\bar{y}^{\delta}-F(\bar{x}_{k_*}^{\delta})\|\le \tau \delta<\|\bar{y}^{\delta}-F(\bar{x}_k^{\delta})\|,~0\le k\le k_{*},
    \end{equation*}
    where $\tau>1$ is a constant. In general, we set $\tau=1.01$.
    Additionally, in Lemma 2.4 of reference \cite{kirsch2011introduction}, we obtain that the Landweber iteration step is the steepest descent step with a stepsize $\zeta_k$ of a quadratic functional $\Psi(\bar{x}^{\delta})=\frac{1}{2}\|\bar{y}^{\delta}-F(\bar{x}^{\delta})\|_2^2$. We denote the descent direction by $\nabla_{\bar{x}^{\delta}}\Psi(\bar{x}^{\delta})$.
    \item \textbf{Method 2}: An improved Landweber iteration regularization. Following the principles of our algorithm, we modify the descent direction using $\nabla_{\bar{x}^{\delta}}\Psi_{\rho}(\bar{x}^{\delta})$, where $\Psi_{\rho}(\bar{x}^{\delta})=\frac{1}{2}\|\bar{y}^{\delta}-F(\bar{x}^{\delta})\|_{\rho}^2$, with $\Vert \cdot \Vert_\rho$ defined in (\ref{de_rho}).
    \item \textbf{Method 3}: Tikhonov regularization method. Our objective is modified to solve the following optimization problem: 
    \begin{equation}
        \bar{x}^{\delta}:=\arg\min\limits_{\bar{x}}\| \bar{y}^\delta-F(\bar{x})\|_2^2+\alpha\|\bar{x}\|_2^2.
    \end{equation}
    Here we set $\alpha=0.005$ or $\alpha=0.1$ in different experiments, and we find that this choice yields the best inversion results for the regularization method. 
\end{enumerate}

Unless otherwise specified, in all numerical examples, for the three traditional regularization methods mentioned above, we fix the maximum number of iterations at 2000, and set the step size $\zeta_k$ in iterations equal to our method ODE-DPS.

\subsection{Inverse problems of heat equation}
\subsubsection{Problem settings}
In this part, we consider the following heat equation:
\begin{equation}\label{5.0}
    \left\{\begin{aligned} 
        u_t(x,y,t) &= a\Delta u(x,y,t) + f(x,y),~(x,y,t)\in\Omega\times (0,T), \\  
        u(x,y,t) &= \psi(x,y,t),~(x,y,t)\in\partial\Omega\times (0,T),\\
        u(x,y,0) &= \phi(x,y),~(x,y)\in\bar{\Omega},
    \end{aligned}\right.
\end{equation}
where $\Omega = [0,1]\times [0,1]$, the coefficient $a$ is a positive constant. 
Next, we will address two inverse problems: the inverse source problem and the inverse initial value problem. Both of these problems are linear and ill-posed, necessitating the use of regularization methods for their solution. We will employ our proposed ODE-DPS inversion algorithm along with three traditional regularization inversion algorithms to tackle these two inverse problems. 

In this example, we focus on the inverse source problem and inverse initial
value problem by measurements
\begin{equation*}
    \bar{y}(x,y) = u(x,y,T).
\end{equation*}
The noisy measurements are generated by adding a random perturbation, i.e.,
\begin{equation}\label{5.1}
\bar{y}^{\delta}:=u^{\delta}(x,y,T)=u(x,y,T)+\varepsilon \eta\max|u(x,y,T)|,
\end{equation}
where $\eta\sim\mathcal{N}(0,I)$, $\varepsilon$ represent the noise percentage.
The corresponding noise level is calculated by $\delta=\|\bar{y}^{\delta}-\bar{y}\|$.

\subsubsection{Implementation details}
To solve the inverse problem, it is necessary to solve the forward problem multiple times. For equation (\ref{5.0}), we utilize a five-point implicit difference scheme. In the inverse source problem, we set $T=1$, the time grid step size as $\Delta t=0.01$, the spatial grid step size as $\Delta x=\Delta y=\frac{1}{63}$, and the regularization parameter of Tikhonov regularization meth as $\alpha=0.005$. In the inverse initial value problem, we set $T = 0.01$, $\Delta t=0.0001$, $\Delta x=\Delta y=\frac{1}{63}$, $\alpha=0.1$. For the noisy measurement data, we set $\varepsilon=0.05$. 

Before applying our method for solving the inverse problems, we trained a new diffusion model with the method in \cite{nichol2021improved} by new datasets shown bellow.  The training algorithm is provided in \ref{A}. The training dataset is composed of functions like 
\begin{equation*}
    \sum_{k=1}^{n_x}\sum_{m=1}^{n_y}\alpha_{k,m}^1\sin k\pi x\sin m\pi y+ \alpha_{k,m}^2\sin k\pi x\cos m\pi y+\alpha_{k,m}^3\cos k\pi x\sin m\pi y+\alpha_{k,m}^4\cos k\pi x\cos m\pi y,
\end{equation*}
where $\alpha_{k,m}^1, \alpha_{k,m}^2,\alpha_{k,m}^3,\alpha_{k,m}^4$ are uniformly random samples in $[-\frac{1}{k^2m^2}, \frac{1}{k^2m^2}]$, $n_x=n_y=15$.
In all experiments of this paper, we select 5000 samples to train the neural network of the score-based generative model.

\subsubsection{Inverse heat source problem}
\begin{figure}[t]
    \centering
    \subfloat[The true source $f_1$(left) and measurements data(right).]{
        \includegraphics[width=0.24\columnwidth,height=.21\linewidth]{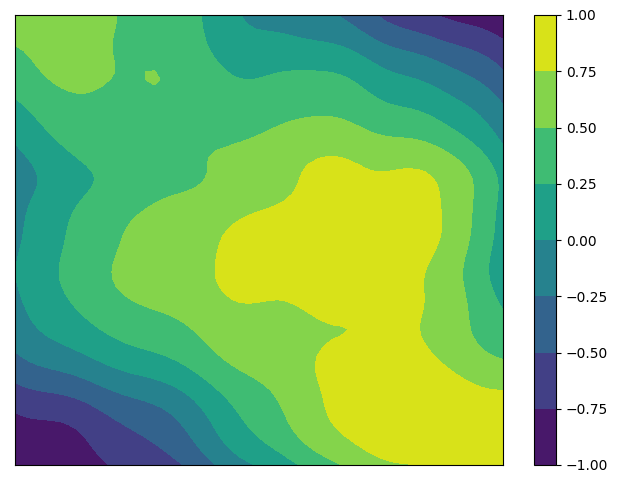}
        \includegraphics[width=.24\columnwidth,height=.21\linewidth]{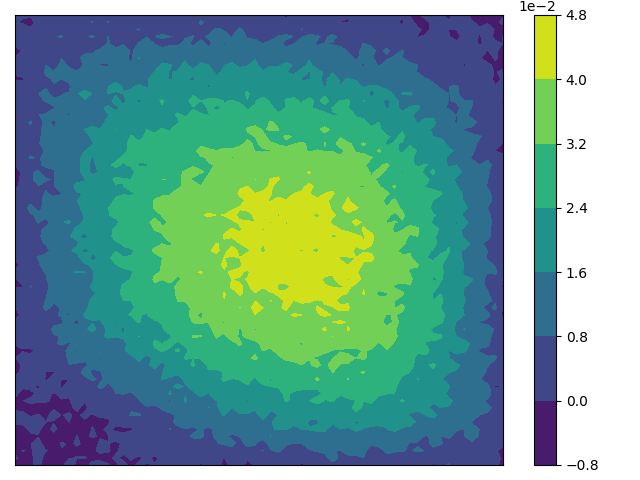}
    }
    \subfloat[The inversion $f_1$ (left) and the errors(right) by ODE-DPS.]{
        \includegraphics[width=0.24\columnwidth,height=.21\linewidth]{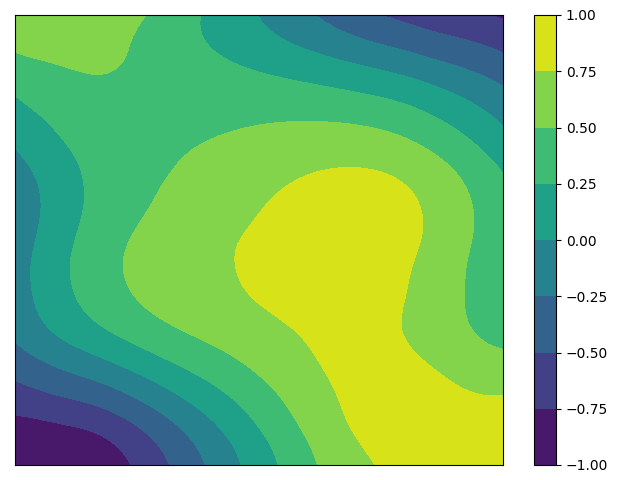}
        \includegraphics[width=.24\columnwidth,height=.21\linewidth]{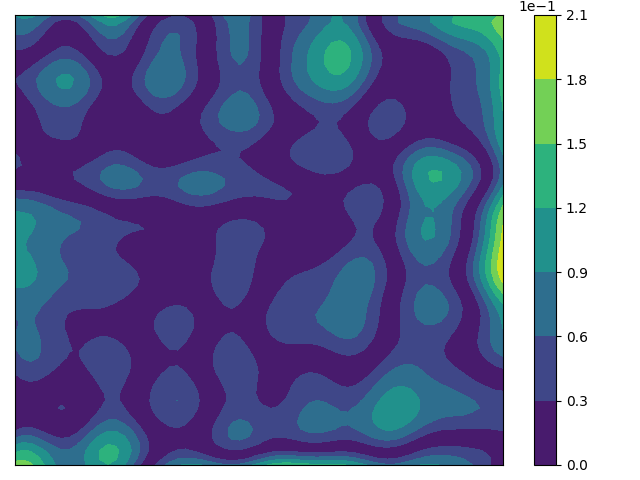}
    }
    
    \subfloat[The inversion $f_1$ (left) and the errors(right) by Method 1(Landweber iteration).]{
        \includegraphics[width=0.24\columnwidth,height=.21\linewidth]{picture/f-gd-sum/f_pre.png}
        \includegraphics[width=.24\columnwidth,height=.21\linewidth]{picture/f-gd-sum/diff_f.png}
    }
    \subfloat[The inversion $f_1$ (left) and the errors(right) by Method 2(Improved Landweber iteration).]{
        \includegraphics[width=0.24\columnwidth,height=.21\linewidth]{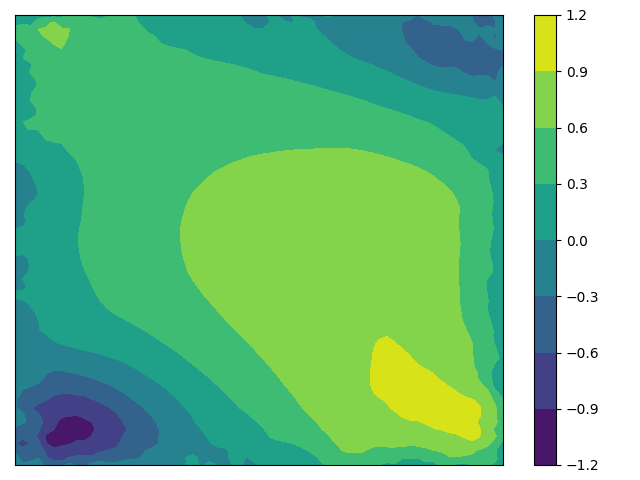}
        \includegraphics[width=.24\columnwidth,height=.21\linewidth]{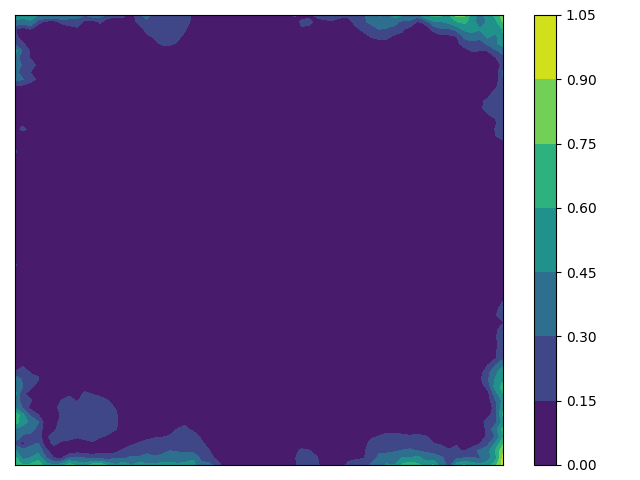}
    }

    \subfloat[The inversion $f_1$ (left) and the errors(right) by Method 3(Tikhonov regularization).]{
        \includegraphics[width=0.24\columnwidth,height=.21\linewidth]{picture/f-gda-sum/f_pre.png}
        \includegraphics[width=.24\columnwidth,height=.21\linewidth]{picture/f-gda-sum/diff_f.png}
    }
    \subfloat[The $e_{f}=\log{\frac{\|f_1-f_{1,rec}\|_2}{\|f_1\|_2}}$.]{\label{f1_f}
        \includegraphics[width=0.24\columnwidth,height=.21\linewidth]{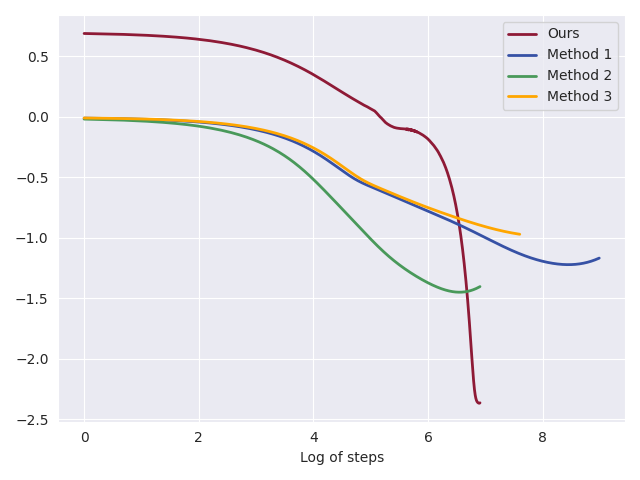}
    }
    \subfloat[The $e_{u}=\log{\frac{\|u^{\delta}(x,y,T;f_1)-u(x,y,T;f_{1,rec})\|_2}{\|u^{\delta}(x,y,T;f_1)\|_2}}$.]{\label{f1_g}
        \includegraphics[width=0.24\columnwidth,height=.21\linewidth]{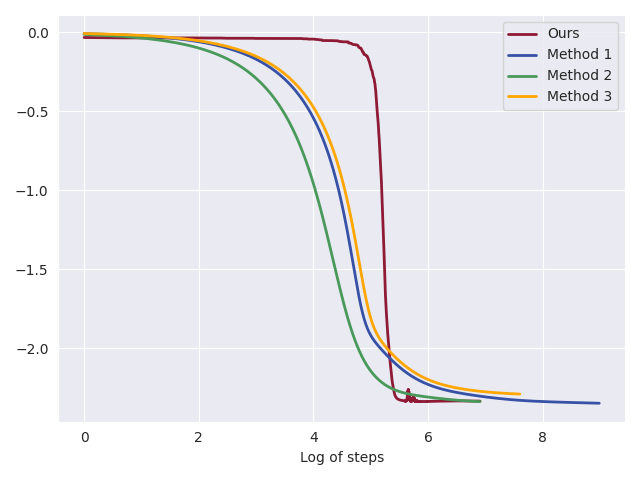}
    }

    \caption{The inversion results of source $f_1$. (a): the true $f_1$ and measurement data $u^{\delta}(x,y,T;f_1)$. (b)(c)(d)(e): the inversion $f_1$ and errors using different regularization methods. (f): the errors between the true $f_1$ and the inversion $f_{1,rec}$. (g): the errors between the measurement data $u^{\delta}(x,y,T;f_1)$ and $u(x,y,T;f_{1,rec})$. Here, we set the noise level $\varepsilon=0.05$ in (\ref{5.1}), the stopping parameter $\tau=1.01$ in discrepancy principle, the regularization parameter $\alpha=0.005$ in Tikhonov regularization.}
    \label{f1_figure}
\end{figure}

\begin{figure}[htbp]
    \centering
    \subfloat[The true source $f_2$(left) and measurements data(right).]{
        \includegraphics[width=0.24\columnwidth,height=.21\linewidth]{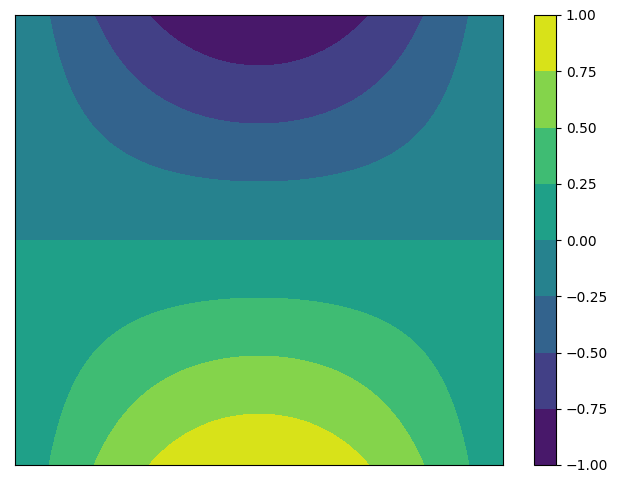}
        \includegraphics[width=.24\columnwidth,height=.21\linewidth]{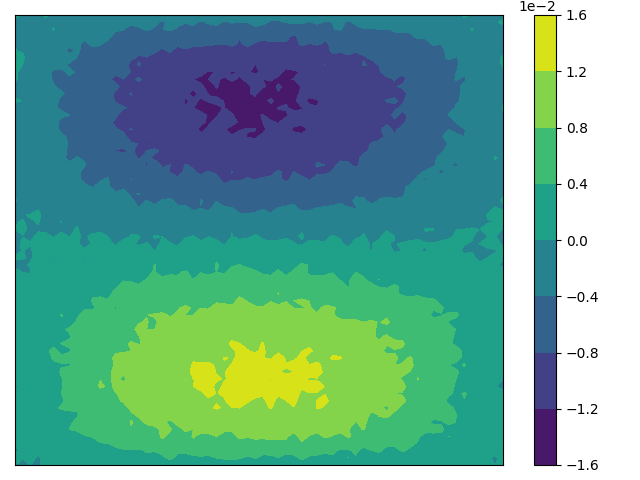}
    }
    \subfloat[The inversion $f_2$ (left) and the errors(right) by ODE-DPS.]{
        \includegraphics[width=0.24\columnwidth,height=.21\linewidth]{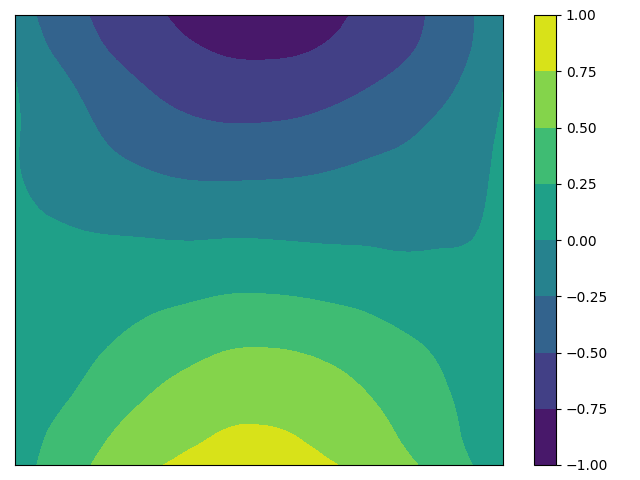}
        \includegraphics[width=.24\columnwidth,height=.21\linewidth]{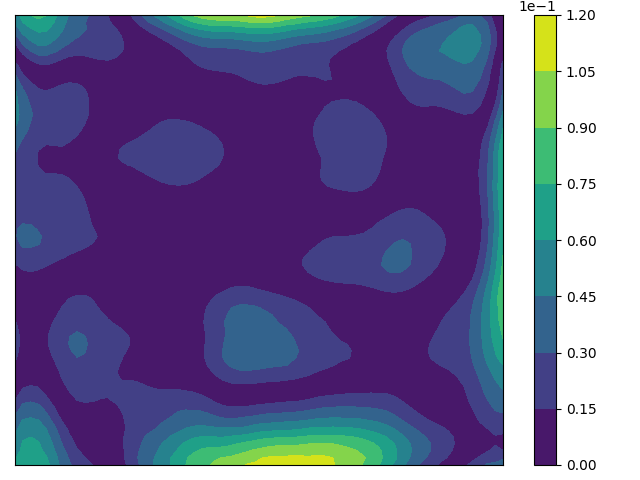}
    }
    
    \subfloat[The inversion $f_2$ (left) and the errors(right) by Method 1(Landweber iteration).]{
        \includegraphics[width=0.24\columnwidth,height=.21\linewidth]{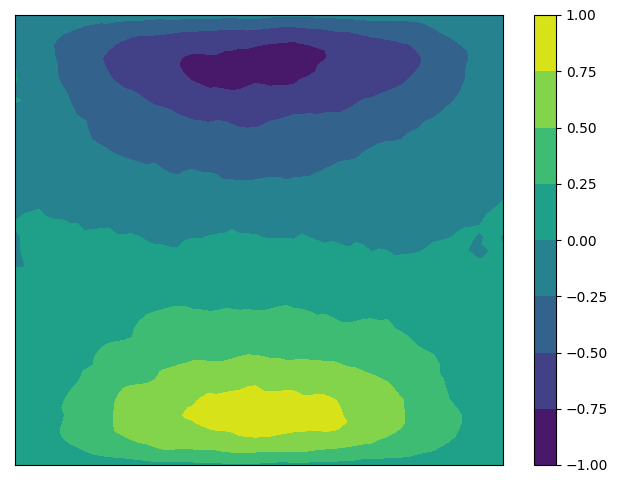}
        \includegraphics[width=.24\columnwidth,height=.21\linewidth]{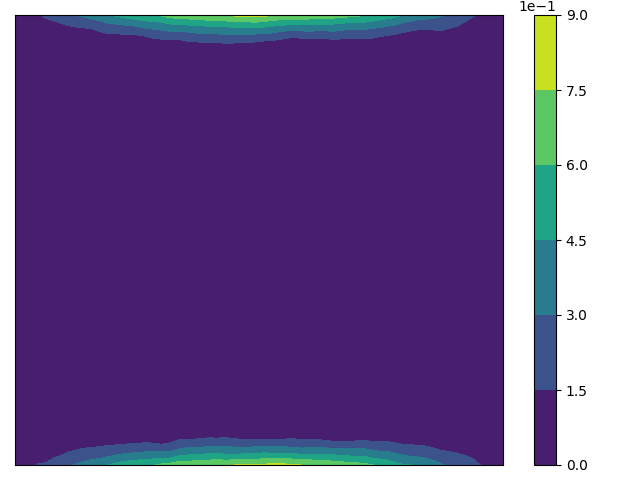}
    }
    \subfloat[The inversion $f_2$ (left) and the errors(right) by Method 2(Improved Landweber iteration).]{
        \includegraphics[width=0.24\columnwidth,height=.21\linewidth]{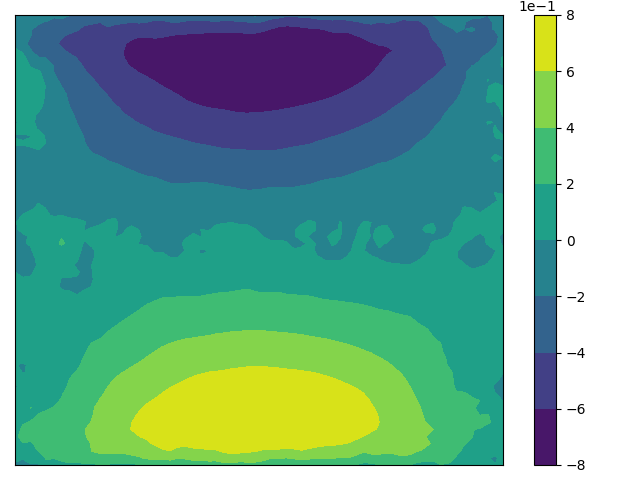}
        \includegraphics[width=.24\columnwidth,height=.21\linewidth]{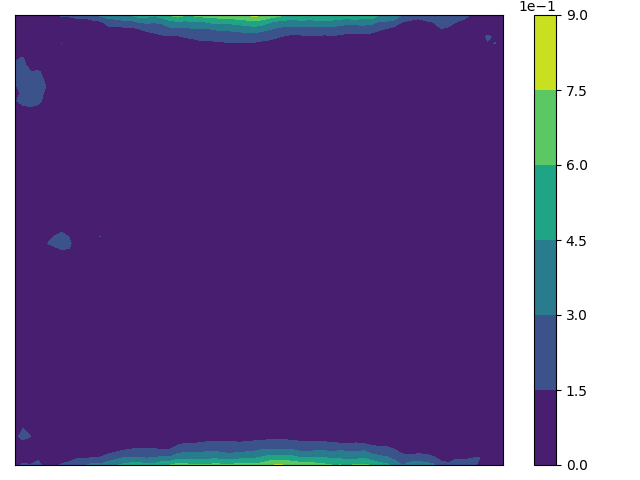}
    }

    \subfloat[The inversion $f_2$ (left) and the errors(right) by Method 3(Tikhonov regularization).]{
        \includegraphics[width=0.24\columnwidth,height=.21\linewidth]{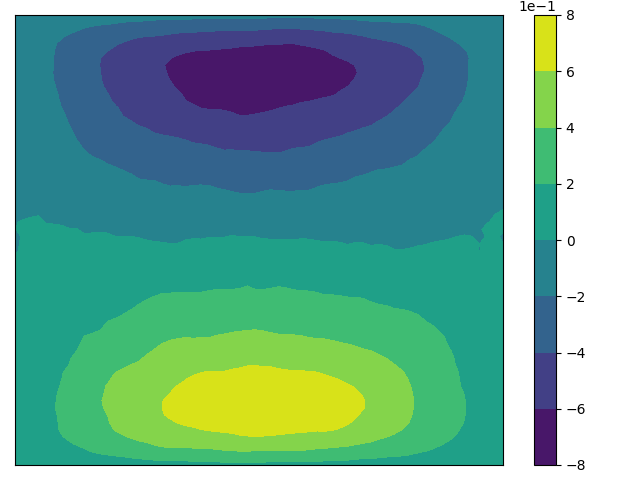}
        \includegraphics[width=.24\columnwidth,height=.21\linewidth]{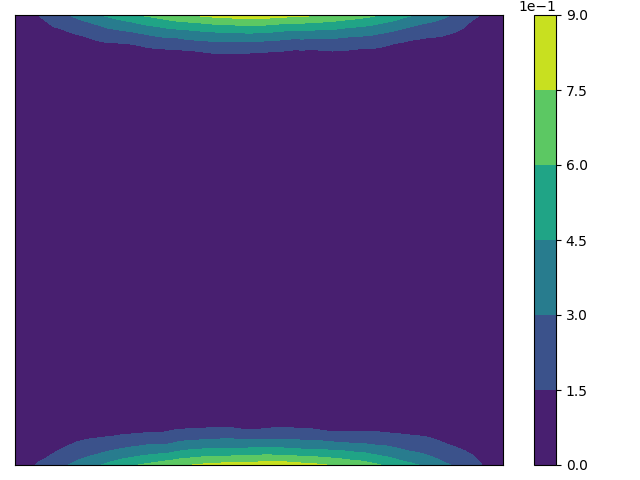}
    }
    \subfloat[The $e_{f}=\log{\frac{\|f_2-f_{2,rec}\|_2}{\|f_2\|_2}}$.]{\label{f2_f}
        \includegraphics[width=0.24\columnwidth,height=0.21\linewidth]{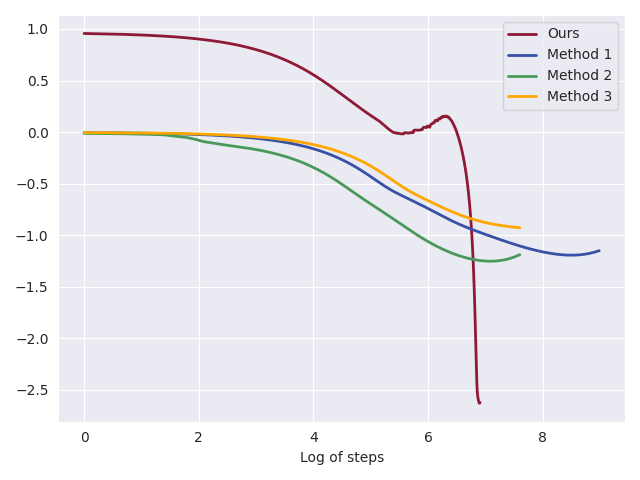}
    }
    \subfloat[The $e_{u}=\log{\frac{\|u^{\delta}(x,y,T;f_2)-u(x,y,T;f_{2,rec})\|_2}{\|u^{\delta}(x,y,T;f_2)\|_2}}$.]{\label{f2_g}
        \includegraphics[width=0.24\columnwidth,height=0.21\linewidth]{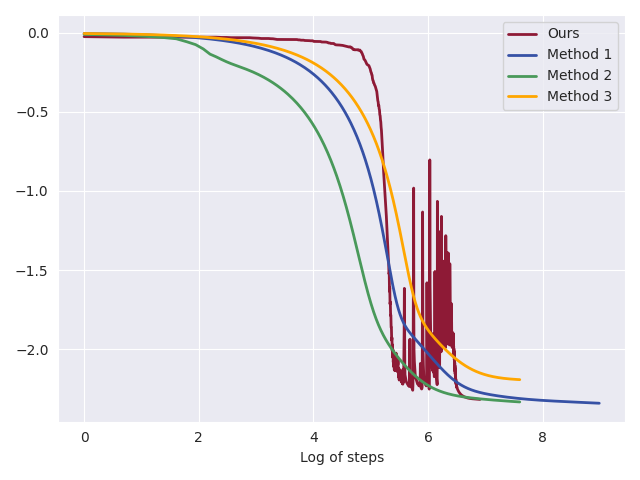}
    }
     \caption{The inversion results of source $f_2$. (a): the true $f_2$ and measurement data $u^{\delta}(x,y,T;f_2)$. (b)(c)(d)(e): the inversion $f_2$ and errors using different regularization methods. (f): the errors between the true $f_2$ and the inversion $f_{2,rec}$. (g): the errors between the measurement data $u^{\delta}(x,y,T;f_2)$ and $u(x,y,T;f_{2,rec})$. Here, we set the noise level $\varepsilon=0.05$ in (\ref{5.1}), the stopping parameter $\tau=1.01$ in discrepancy principle, the regularization parameter $\alpha=0.005$ in Tikhonov regularization.}
    \label{f2_figure}
\end{figure}

\begin{figure}[t]
    \centering
    \subfloat[The true source $f_3$(left) and measurements data(right).]{
        \includegraphics[width=0.24\columnwidth,height=.21\linewidth]{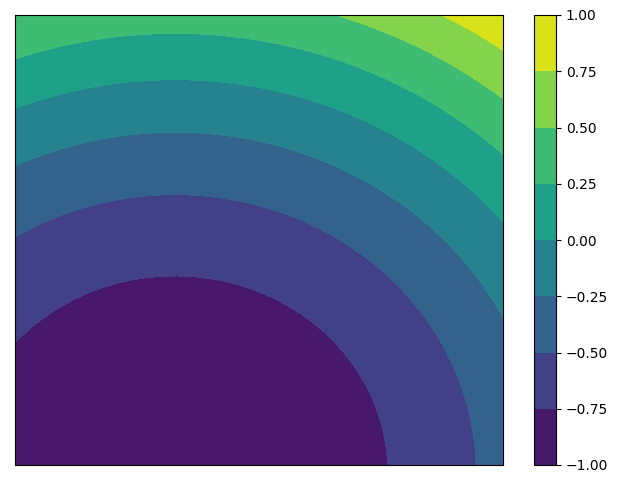}
        \includegraphics[width=.24\columnwidth,height=.21\linewidth]{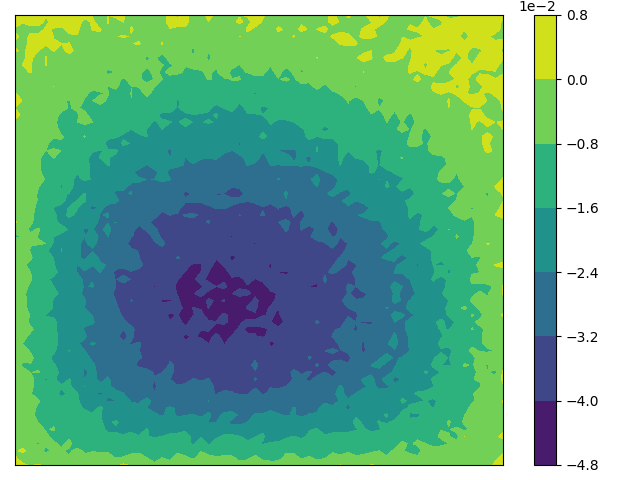}
    }
    \subfloat[The inversion $f_3$ (left) and the errors(right) by ODE-DPS.]{
        \includegraphics[width=0.24\columnwidth,height=.21\linewidth]{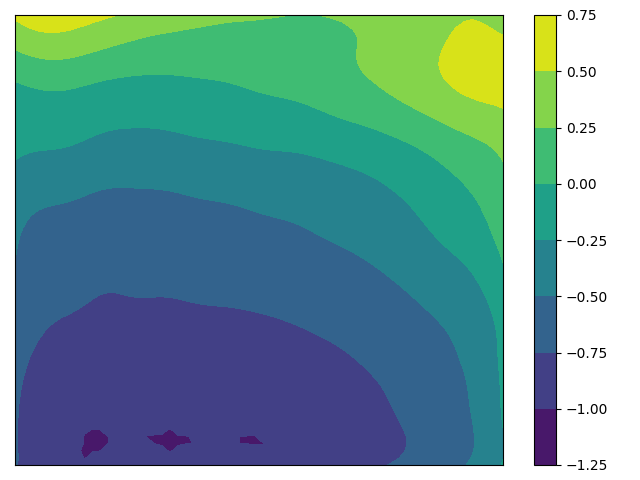}
        \includegraphics[width=.24\columnwidth,height=.21\linewidth]{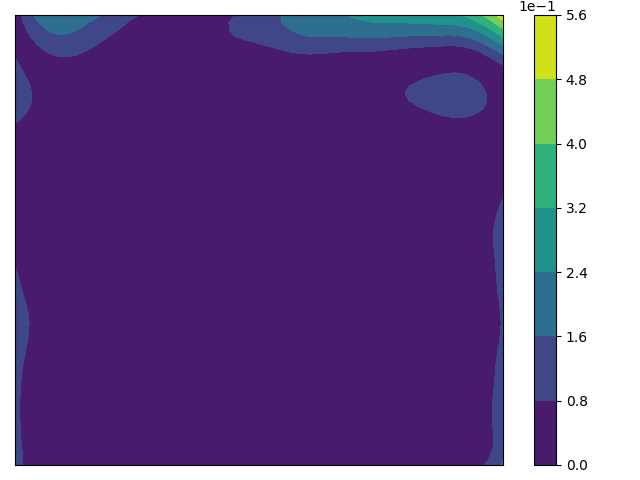}
    }
    
    \subfloat[The inversion $f_3$ (left) and the errors(right) by Method 1(Landweber iteration).]{
        \includegraphics[width=0.24\columnwidth,height=.21\linewidth]{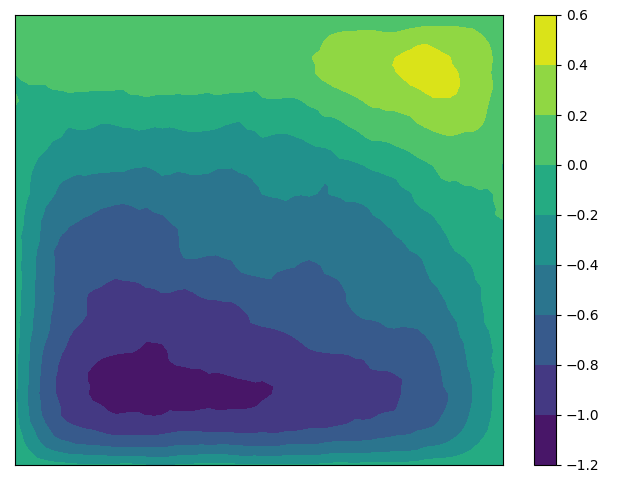}
        \includegraphics[width=.24\columnwidth,height=.21\linewidth]{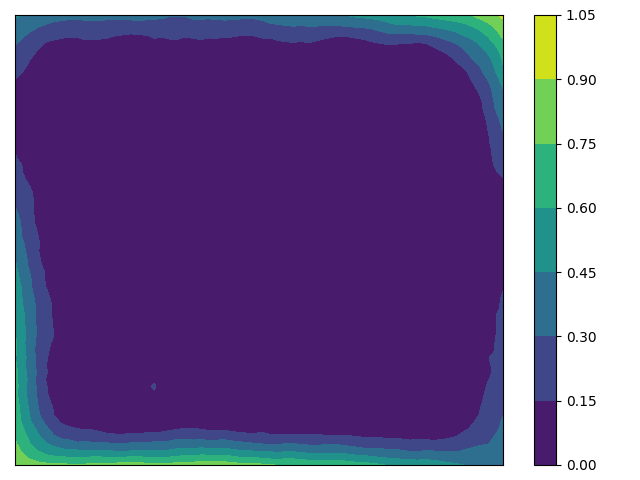}
    }
    \subfloat[The inversion $f_3$ (left) and the errors(right) by Method 2(Improved Landweber iteration).]{
        \includegraphics[width=0.24\columnwidth,height=.21\linewidth]{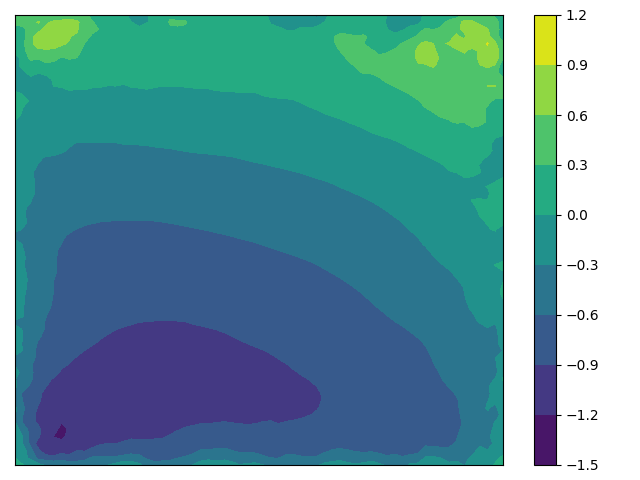}
        \includegraphics[width=.24\columnwidth,height=.21\linewidth]{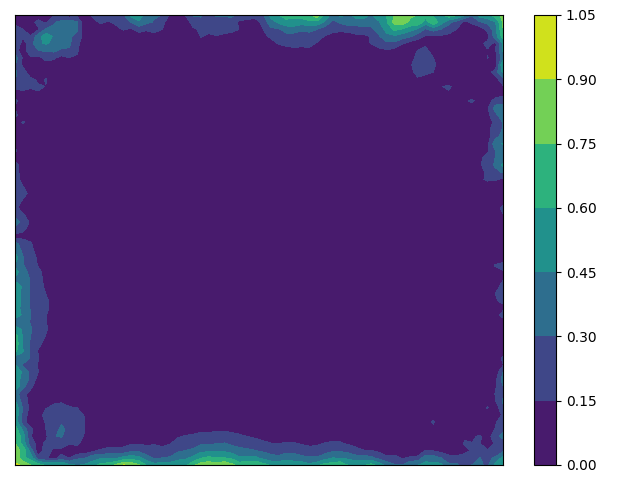}
    }

    \subfloat[The inversion $f_3$ (left) and the errors(right) by Method 3(Tikhonov regularization).]{
        \includegraphics[width=0.24\columnwidth,height=.21\linewidth]{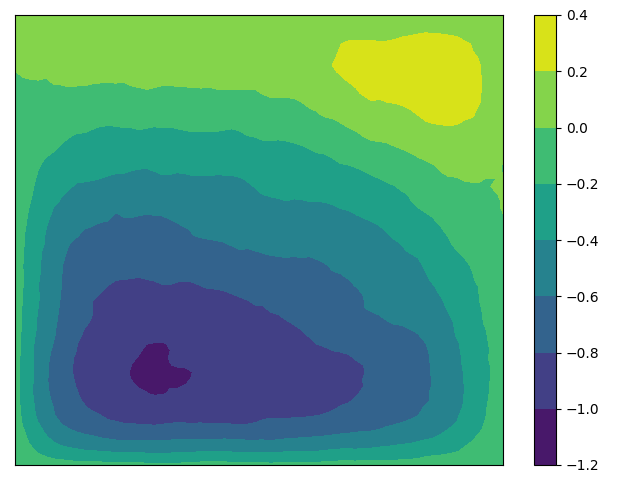}
        \includegraphics[width=.24\columnwidth,height=.21\linewidth]{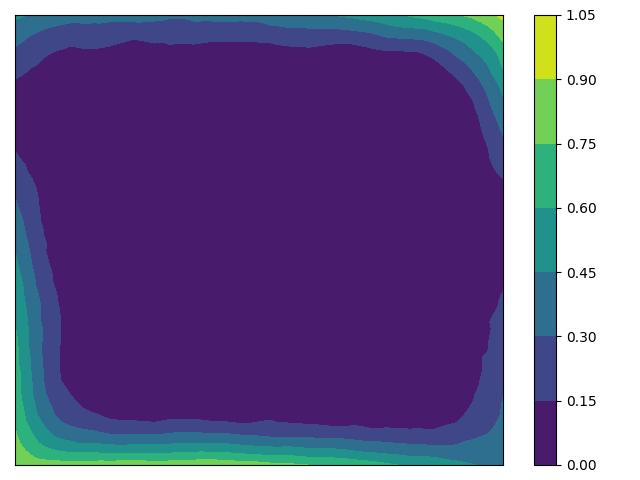}
    }
    \subfloat[The $e_{f}=\log{\frac{\|f_3-f_{3,rec}\|_2}{\|f_3\|_2}}$.]{\label{f3_f}
        \includegraphics[width=0.24\columnwidth,height=.21\linewidth]{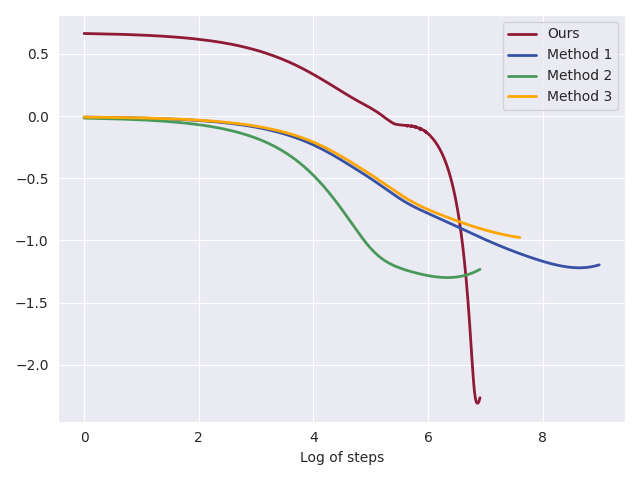}
    }
    \subfloat[The $e_{u}=\log{\frac{\|u^{\delta}(x,y,T;f_3)-u(x,y,T;f_{3,rec})\|_2}{\|u^{\delta}(x,y,T;f_3)\|_2}}$.]{\label{f3_g}
        \includegraphics[width=0.24\columnwidth,height=.21\linewidth]{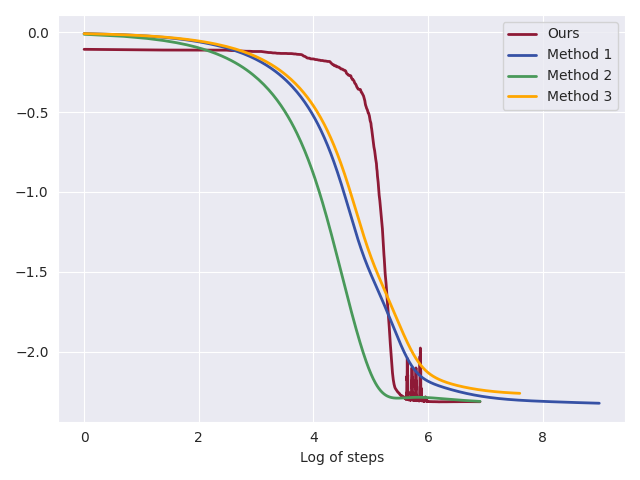}
    }
    \caption{The inversion results of source $f_3$. (a): the true $f_3$ and measurement data $u^{\delta}(x,y,T;f_3)$. (b)(c)(d)(e): the inversion $f_3$ and errors using different regularization methods. (f): the errors between the true $f_3$ and the inversion $f_{3,rec}$. (g): the errors between the measurement data $u^{\delta}(x,y,T;f_3)$ and $u(x,y,T;f_{3,rec})$. Here, we set the noise level $\varepsilon=0.05$ in (\ref{5.1}), the stopping parameter $\tau=1.01$ in discrepancy principle, the regularization parameter $\alpha=0.005$ in Tikhonov regularization.}
    \label{f3_figure}
\end{figure}
In this example, we take the functions of problem (\ref{5.0}) as follows: the diffusion coefficient $a=1$, the boundary $\psi(x,y,t)=0$ for $(x,y,t)\in\partial\Omega\times (0,T)$, the initial value $\phi(x,y)=20\sin 2\pi x\sin 2\pi y$ for $(x,y)\in\Omega$.
The parameters of our algorithm ODE-DPS are set as $\gamma=1,~ \zeta=9,~ c=1e-8, ~N=1000$. 

Next, we will explore the inversion of three distinct source terms to evaluate the performance of our inversion technique. To ascertain the effectiveness of our proposed method, we will select one source term that mirrors those present in the training dataset but is absent from it. Furthermore, we will opt for two additional source functions that exhibit significant deviation from those featured in the training dataset. In the following, we will consider three representative examples given by:
\begin{itemize}
  \item \textbf{Case 1:} We select source $f_1$ shared the same form of training dataset, but $f_1$ is not included in the training set. The ground truth source $f_1$ is shown in Figure \ref{f1_figure}.
  \item \textbf{Case 2:} We choose $f_2(x,y)=(1-2x)\sin\pi y$, which represents a polynomial function multiplied by a trigonometric function. The ground truth source $f_2$ is shown in Figure \ref{f2_figure}.
  \item \textbf{Case 3:} We take $f_3(x,y)=\frac{7}{5}(x^2+(y-\frac{1}{3})^2)-1$. The ground truth source $f_3$ is shown in Figure \ref{f3_figure}.
\end{itemize}

\textit{\textbf{Results and discussion:}} 
We compare the performance of the ODE-DPS regularization inversion method with that of the Landweber iteration regularization (Method 1), an improved Landweber iteration regularization (Method 2) and Tikhonov regularization (Method 3). All hyperparameters in the tested methods are manually tuned for optimal performance. Table \ref{heat_source_table} illustrates the relative $l_2$ error for different heat sources $f$ using various regularization methods. Visual comparisons can be found in Figures \ref{f1_figure}-\ref{f3_figure}. We summarize our findings as follows:
\begin{enumerate}
    \item For all three numerical examples, the ODE-DPS inversion algorithm yielded the best numerical inversion results. Furthermore, our proposed inversion algorithm required fewer iterations compared to the Landweber iterative regularization method (Method 1) and Tikhonov regularization method (Method 3).
    \item The improved Landweber iterative regularization method, referred to as Method 2, produced inversion outcomes with reduced relative $l_2$ error in comparison to the original Landweber iterative regularization method (Method 1). This suggests that adjusting the gradient descent direction using $\nabla_{\bar{x}^{\delta}}\Psi_{\rho}(\bar{x}^{\delta})$ throughout the iterations of the inversion algorithm enhances its performances.
    \item Figures \ref{f1_f}, \ref{f2_f}, and \ref{f3_f} indicate that traditional regularization methods without regularization item lead to semi-convergent numerical inversion results with increasing iterations while Tikhonov regularization method and our proposed algorithm shows a consistent decrease in relative $l_2$ error after a certain number of iterations. Furthermore, Figures \ref{f1_g}, \ref{f2_g}, and \ref{f3_g} demonstrate that all the inversion methods produce source terms $f$ of varying accuracy but their respective numerical solutions $u(x,y,T;f_{i,rec})$ and measurement data $u^{\delta}(x,y,T;f_{i})$ show nearly identical discrepancies, where $i=1,2,3$.
\end{enumerate}

\begin{table}[ht]
    \begin{center}
    \begin{tabular}{llllllr}
    \toprule[1.pt]
     &\multicolumn{2}{c}{$f_1$} & \multicolumn{2}{c}{$f_2$} & \multicolumn{2}{c}{$f_3$} \\
    \cline{2-3} \cline{4-5} \cline{6-7}
    method      & steps & relative $l_2$ error  & steps & relative $l_2$ error & steps & relative $l_2$ error\\ \hline
    Method 1  &  1471 &  34.4\% & 1704 & 34.1\% & 1732 & 33.9\% \\ 
    Method 2     & 523 & 24.1\% & 817 & 29.4\% & 637 & 27.4\%\\ 
    Method 3     &  2000 & 37.9\% & 2000 & 39.5\% & 2000 & 37.7\%\\
    ODE-DPS(ours) & 1000 &\textbf{ 9.4\%} & 1000 & \textbf{7.3\%} & 1000 & \textbf{10.4\%} \\ 
    \bottomrule[1.0pt]
    \end{tabular}
    \end{center}
    \caption{The relative $l_2$ error and iteration steps of inverse heat source using different inversion methods. Methods 1, 2, and 3 correspond to Landweber iteration, improved Landweber iteration, and Tikhonov regularization, respectively. We set the noise level $\varepsilon=0.05$ in (\ref{5.1}), the stopping parameter $\tau=1.01$ in discrepancy principle, the regularization parameter $\alpha=0.005$ in Tikhonov regularization.}
    \label{heat_source_table}
\end{table}
\subsubsection{Inverse initial value problem}

In this experiment, we take the functions of problem (\ref{5.0}) as follows: the diffusion coefficient $a=1$, the heat source $f(x,y)=0$ for $(x,y)\in\Omega$.  The parameters of our algorithm ODE-DPS are set as $\gamma=0.75,~ \zeta=2,~ c=1e-5,~ N=1000$.

\begin{figure}[htbp!]
    \centering
    \subfloat[The true initial value $\phi_1$(left) and measurements data(right).]{
        \includegraphics[width=0.24\columnwidth,height=.21\linewidth]{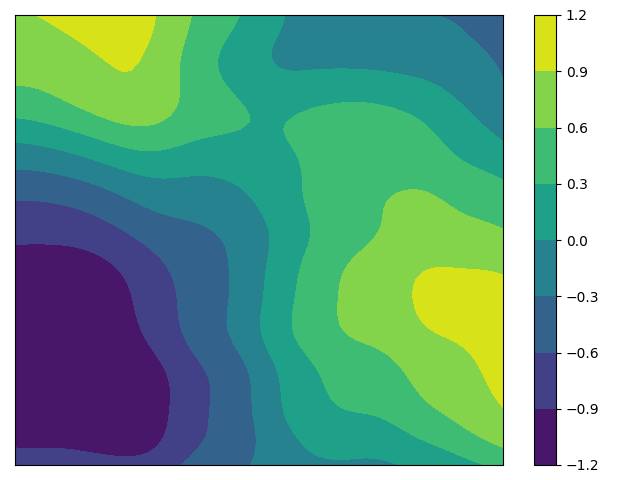}
        \includegraphics[width=.24\columnwidth,height=.21\linewidth]{picture/f-diffusion-sum/label_noise.png}
    }
    \subfloat[The inversion $\phi_1$(left) and errors(right) by ODE-DPS.]{
        \includegraphics[width=0.24\columnwidth,height=.21\linewidth]{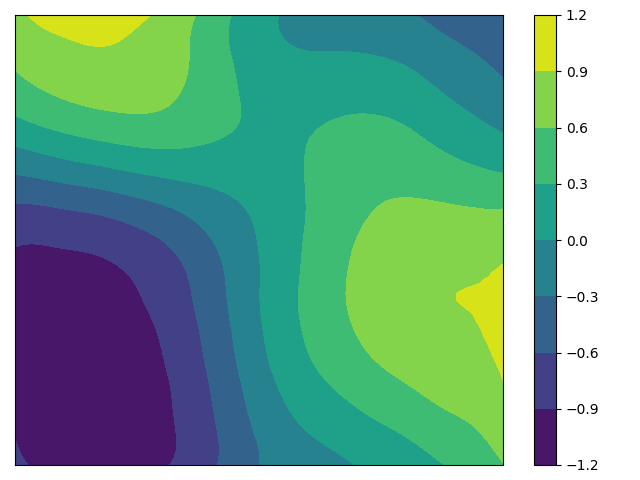}
        \includegraphics[width=.24\columnwidth,height=.21\linewidth]{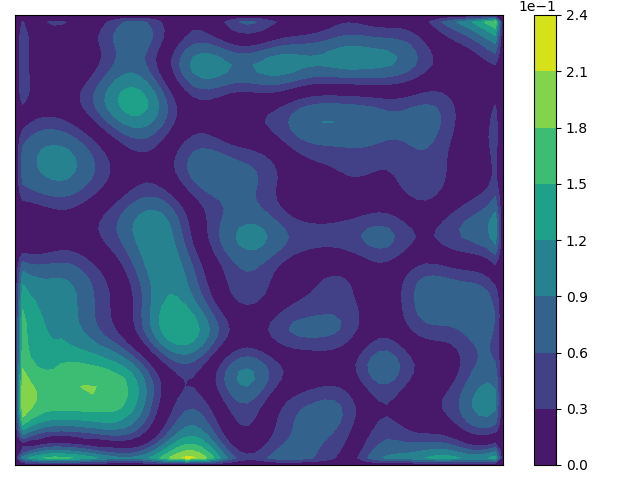}
    }
    
    \subfloat[The inversion $\phi_1$(left) and errors(right) by Method 1(Landweber iteration).]{
        \includegraphics[width=0.24\columnwidth,height=.21\linewidth]{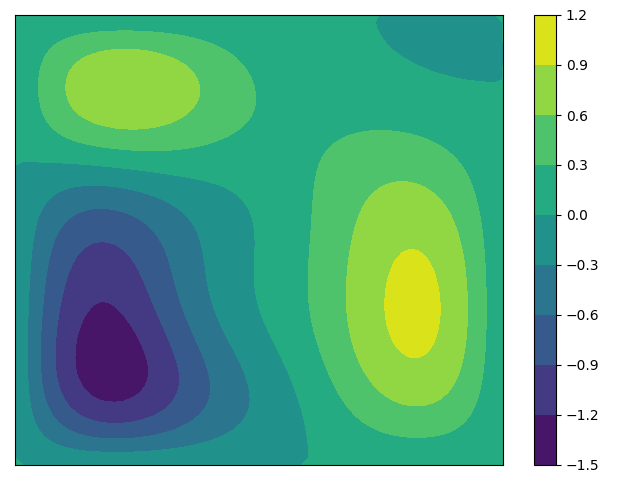}
        \includegraphics[width=.24\columnwidth,height=.21\linewidth]{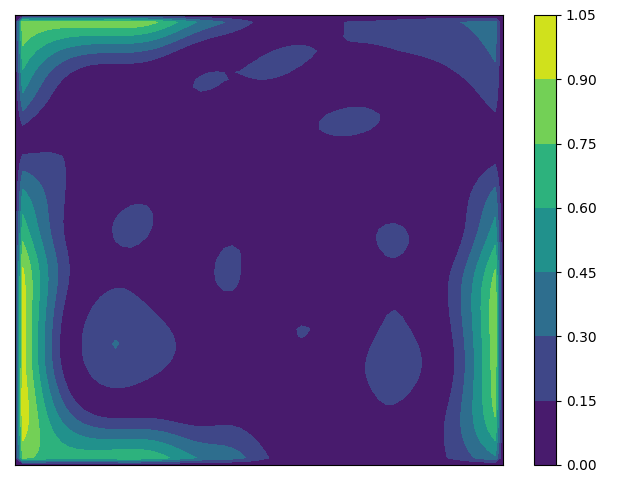}
    }
    \subfloat[The inversion $\phi_1$(left) and errors(right) by Method 2(Improved Landweber iteration).]{
        \includegraphics[width=0.24\columnwidth,height=.21\linewidth]{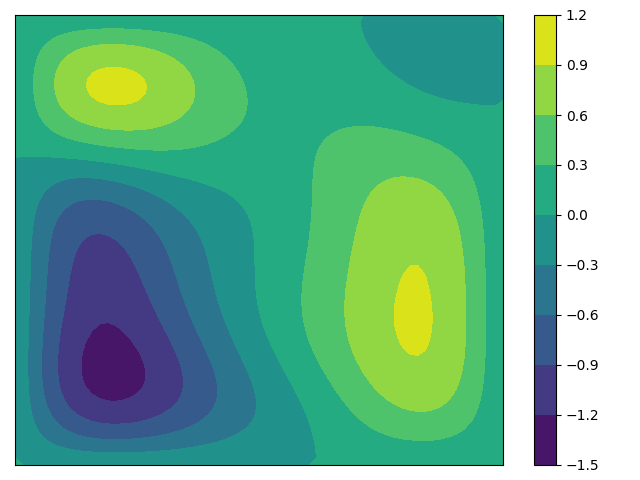}
        \includegraphics[width=.24\columnwidth,height=.21\linewidth]{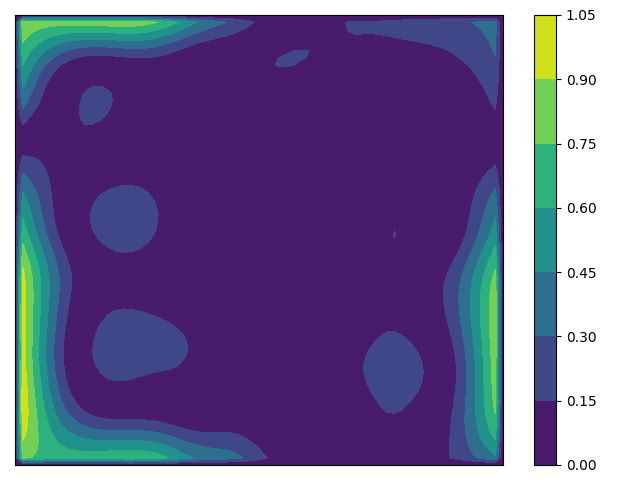}
    }

    \subfloat[The inversion $\phi_1$(left) and errors(right) by Method 3(Tikhonov regularization).]{
        \includegraphics[width=0.24\columnwidth,height=.21\linewidth]{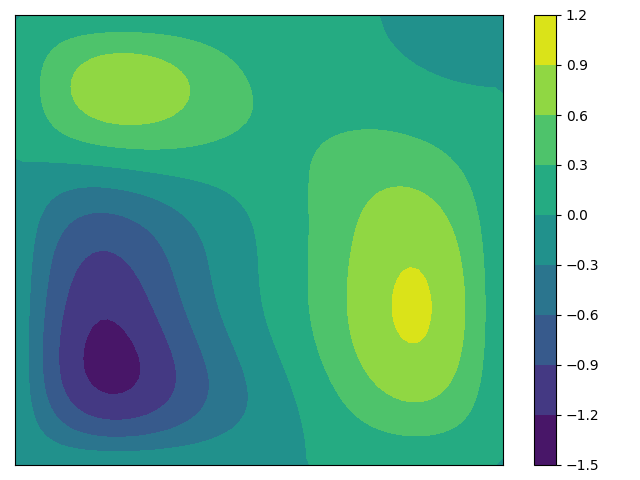}
        \includegraphics[width=.24\columnwidth,height=.21\linewidth]{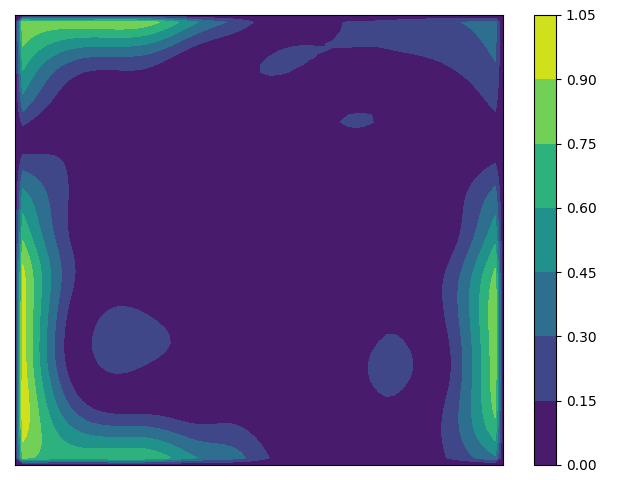}
    }
    \subfloat[The $e_f=\log{\frac{\|\phi_1-\phi_{1,rec}\|_2}{\|\phi_1\|_2}}$.]{
        \label{initial_value_distance_1}
        \includegraphics[width=0.24\columnwidth,height=.21\linewidth]{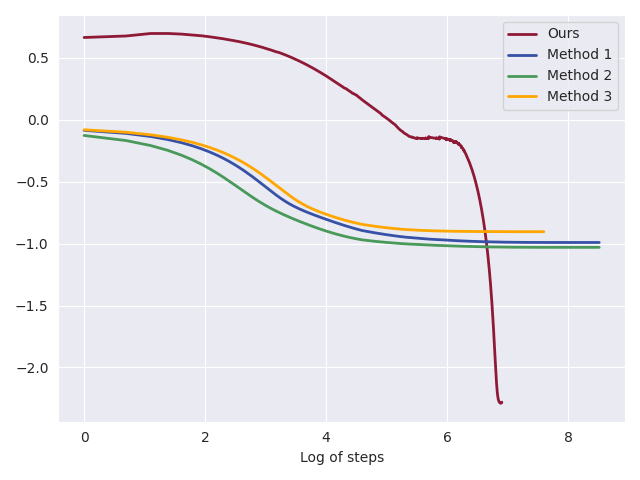}
    }
    \subfloat[The $e_u=\log{\frac{\|u^\delta(x,y,T;\phi_1)-u(x,y,T;\phi_{1,rec})\|_2}{\|u^\delta(x,y,T;\phi_1)\|_2}}$.]{
        \includegraphics[width=0.24\columnwidth,height=.21\linewidth]{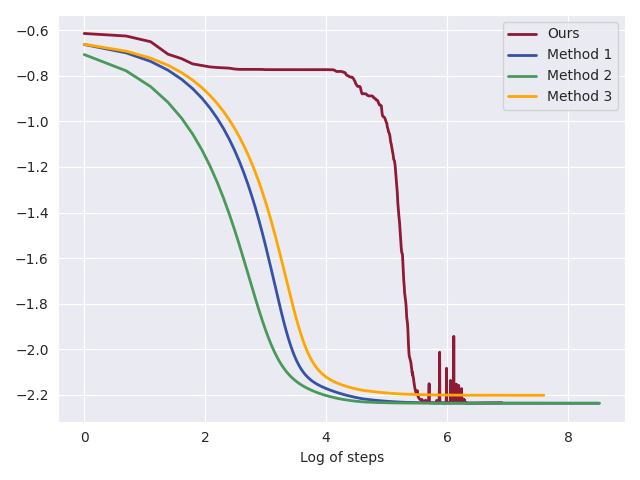}
    }

    \caption{The inversion results of initial value $\phi_1$. (a): the true $\phi_1$ and measurement data $u^\delta(x,y,T; \phi_1)$. (b)(c)(d)(e): the inversion $\phi_1$ and errors using different regularization methods. (f): the errors between the true $\phi_1$ and the inversion $\phi_{1,rec}$. (g): the errors between the measurement data $u^\delta(x,y,T; \phi_1)$ and $u(x,y,T; \phi_{1,rec})$. Here, we set the noise level $\varepsilon=0.05$ in (\ref{5.1}), the stopping parameter $\tau=1.01$ in discrepancy principle, the regularization parameter $\alpha=0.1$ in Tikhonov regularization.}
    \label{phi1_figure}
\end{figure}

\begin{figure}[htbp!]
    \centering
    \subfloat[The true initial value $\phi_2$(left) and measurements data(right).]{
        \includegraphics[width=0.24\columnwidth,height=.21\linewidth]{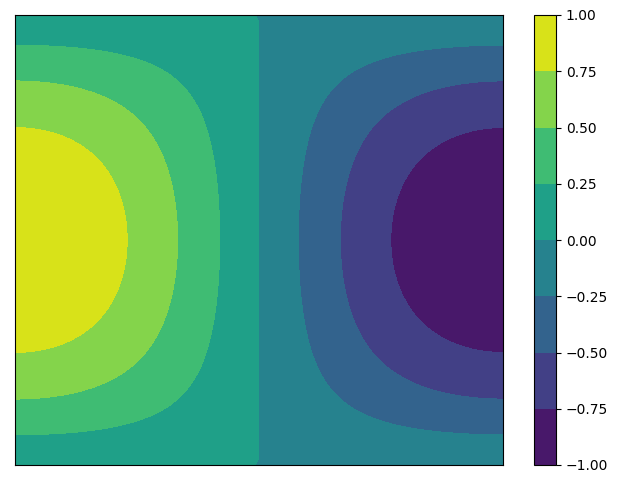}
        \includegraphics[width=.24\columnwidth,height=.21\linewidth]{picture/f-diffusion-continue/label_noise.png}
    }
    \subfloat[The inversion $\phi_2$(left) and errors(right) by ODE-DPS.]{
        \includegraphics[width=0.24\columnwidth,height=.21\linewidth]{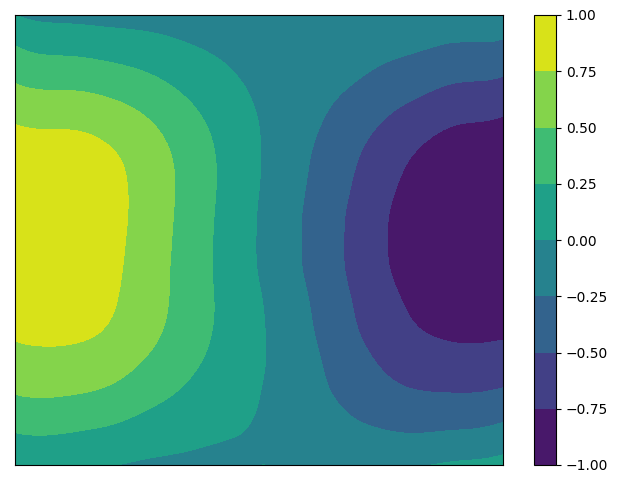}
        \includegraphics[width=.24\columnwidth,height=.21\linewidth]{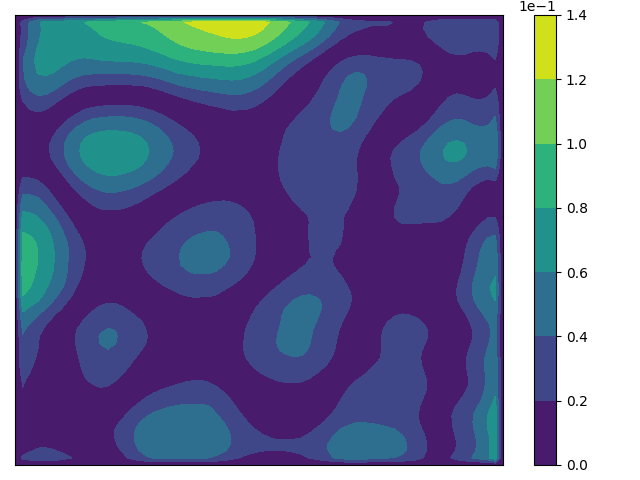}
    }
    
    \subfloat[The inversion $\phi_2$(left) and errors(right) by Method 1(Landweber iteration).]{
        \includegraphics[width=0.24\columnwidth,height=.21\linewidth]{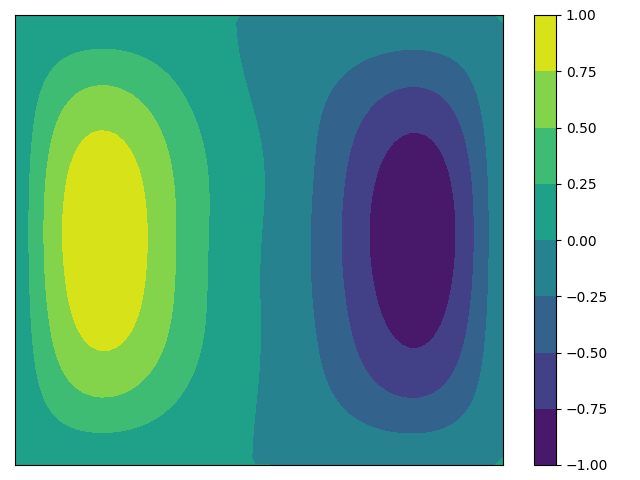}
        \includegraphics[width=.24\columnwidth,height=.21\linewidth]{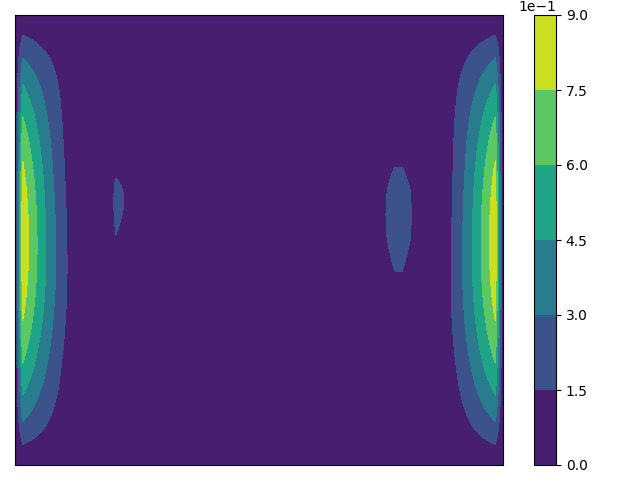}
    }
    \subfloat[The inversion $\phi_2$(left) and errors(right) by Method 2(Improved Landweber iteration).]{
        \includegraphics[width=0.24\columnwidth,height=.21\linewidth]{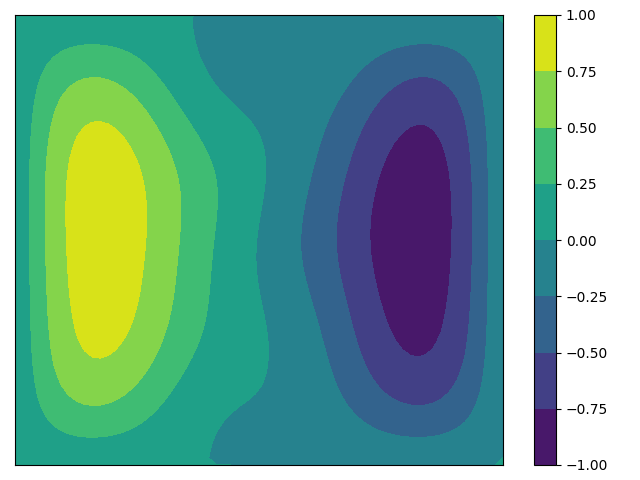}
        \includegraphics[width=.24\columnwidth,height=.21\linewidth]{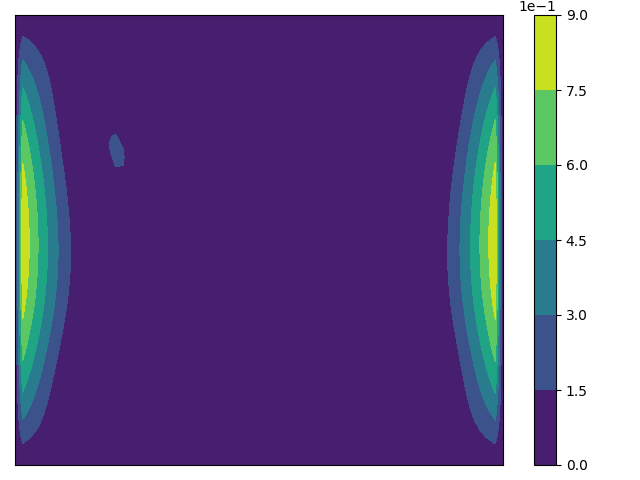}
    }

    \subfloat[The inversion $\phi_2$(left) and errors(right) by Method 3(Tikhonov regularization).]{
        \includegraphics[width=0.24\columnwidth,height=.21\linewidth]{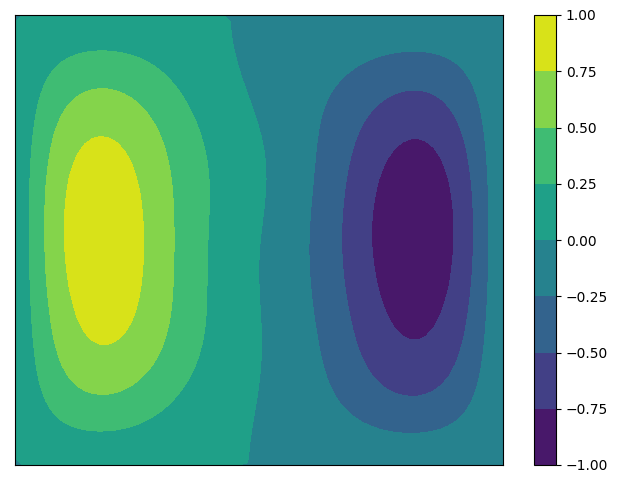}
        \includegraphics[width=.24\columnwidth,height=.21\linewidth]{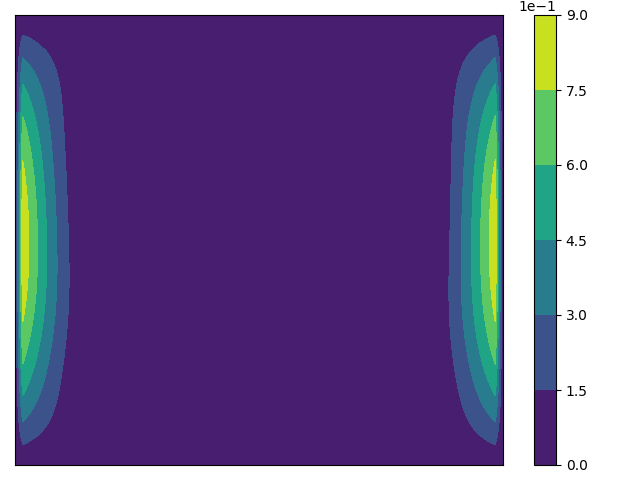}
    }
    \subfloat[The $e_f=\log{\frac{\|\phi_2-\phi_{2,rec}\|_2}{\|\phi_2\|_2}}$.]{
        \label{initial_value_distance_2}
        \includegraphics[width=0.24\columnwidth,height=.21\linewidth]{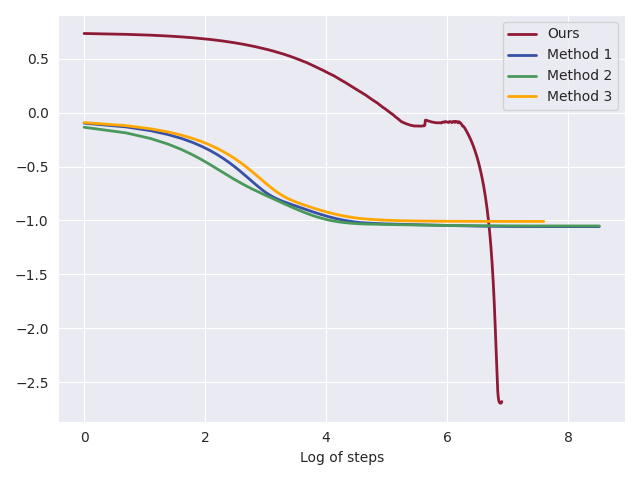}
    }
    \subfloat[The $e_u=\log{\frac{\|u^\delta(x,y,T;\phi_2)-u(x,y,T;\phi_{2,rec})\|_2}{\|u^\delta(x,y,T;\phi_2)\|_2}}$.]{
        \includegraphics[width=0.24\columnwidth,height=.21\linewidth]{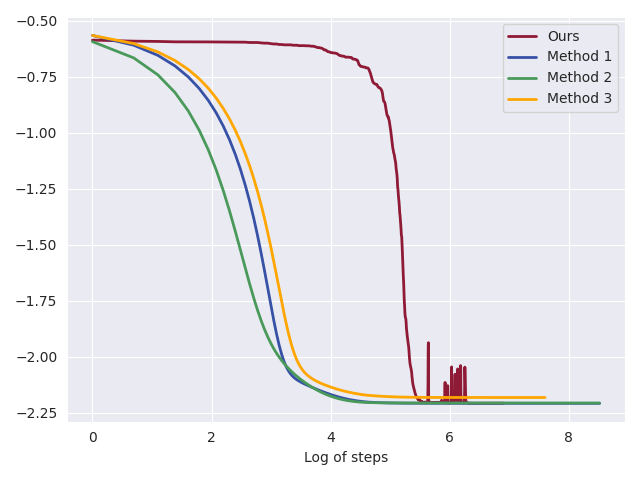}
    }

    \caption{The inversion results of initial value $\phi_2$. (a): the true $\phi_2$ and measurement data $u^\delta(x,y,T; \phi_2)$. (b)(c)(d)(e): the inversion $\phi_2$ and errors using different regularization methods. (f): the errors between the true $\phi_2$ and the inversion $\phi_{2,rec}$. (g): the errors between the measurement data $u^\delta(x,y,T; \phi_2)$ and $u(x,y,T; \phi_{2,rec})$. Here, we set the noise level $\varepsilon=0.05$ in (\ref{5.1}), the stopping parameter $\tau=1.01$ in discrepancy principle, the regularization parameter $\alpha=0.1$ in Tikhonov regularization.}
    \label{phi2_figure}
\end{figure}

\begin{figure}[htbp!]
    \centering
    \subfloat[The true initial value $\phi_3$(left) and measurements data(right).]{
        \includegraphics[width=0.24\columnwidth,height=.21\linewidth]{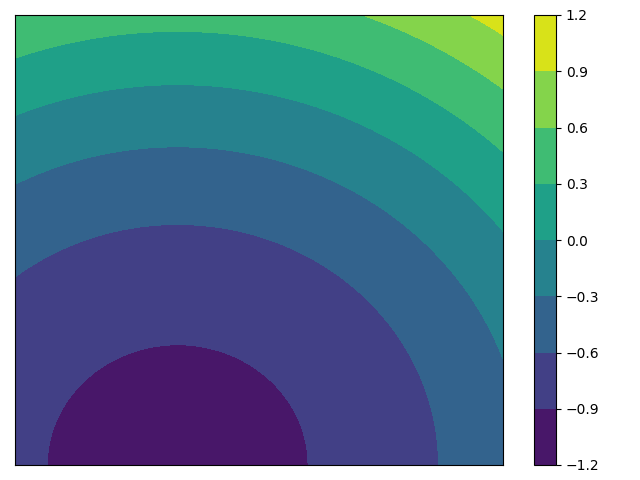}
        \includegraphics[width=.24\columnwidth,height=.21\linewidth]{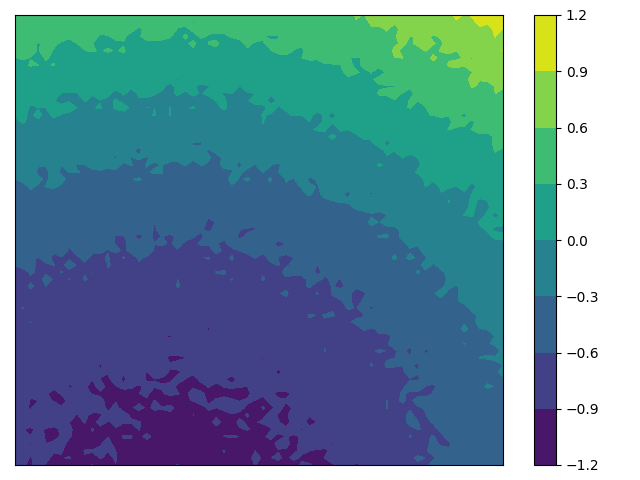}
    }
    \subfloat[The inversion $\phi_3$(left) and errors(right) by ODE-DPS.]{
        \includegraphics[width=0.24\columnwidth,height=.21\linewidth]{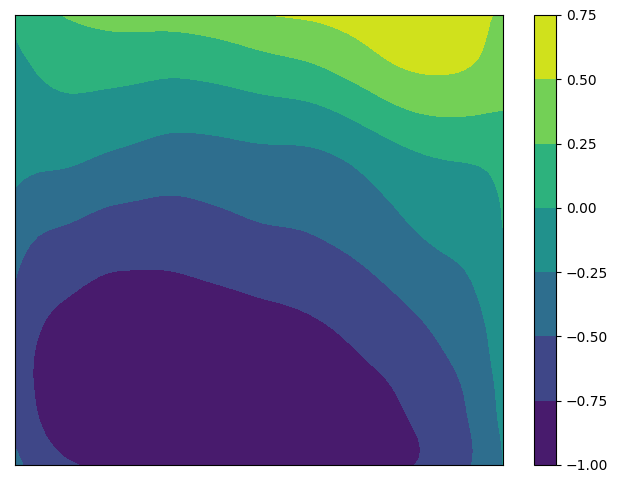}
        \includegraphics[width=.24\columnwidth,height=.21\linewidth]{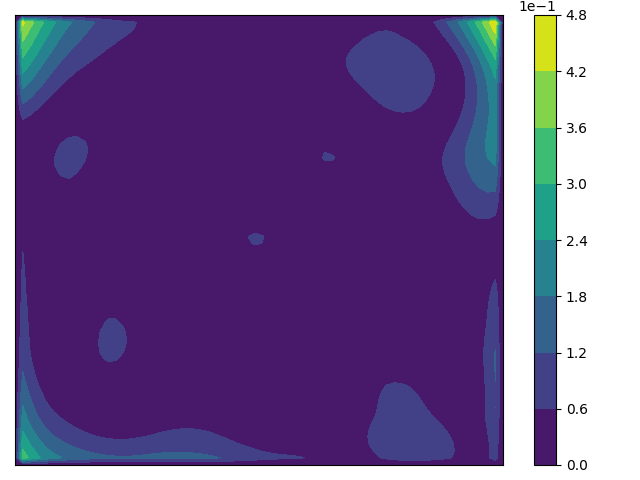}
    }
    
    \subfloat[The inversion $\phi_3$(left) and errors(right) by Method 1(Landweber iteration).]{
        \includegraphics[width=0.24\columnwidth,height=.21\linewidth]{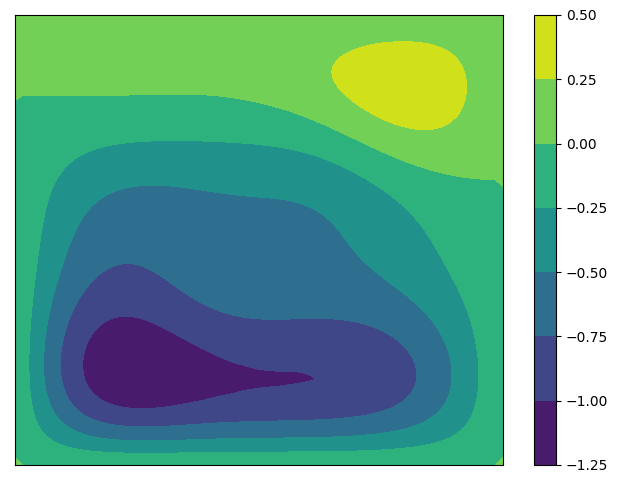}
        \includegraphics[width=.24\columnwidth,height=.21\linewidth]{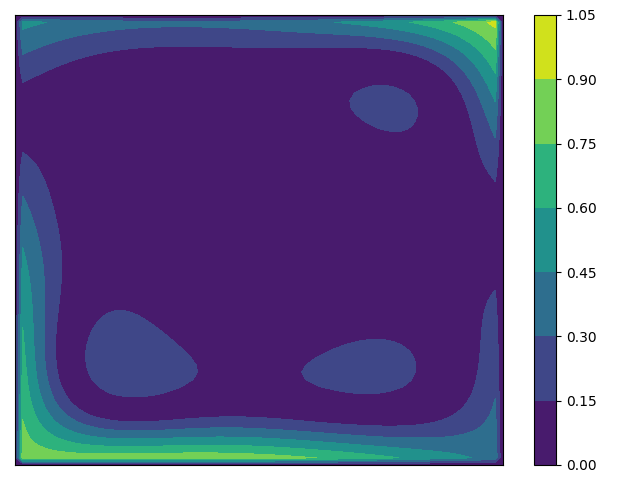}
    }
    \subfloat[The inversion $\phi_3$(left) and errors(right) by Method 2(Improved Landweber iteration).]{
        \includegraphics[width=0.24\columnwidth,height=.21\linewidth]{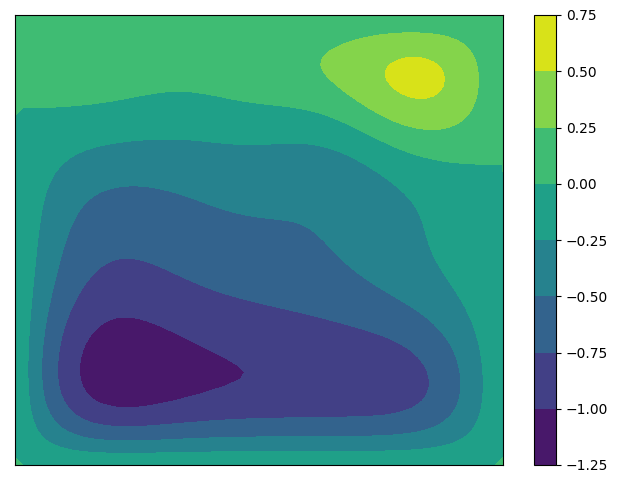}
        \includegraphics[width=.24\columnwidth,height=.21\linewidth]{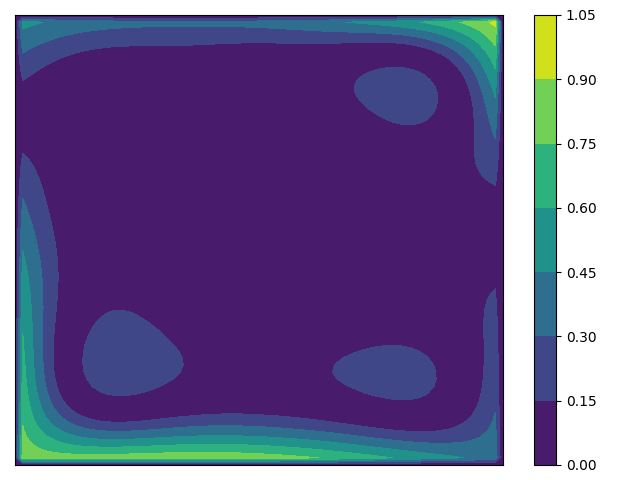}
    }

    \subfloat[The inversion $\phi_3$(left) and errors(right) by Method 3(Tikhonov regularization).]{
        \includegraphics[width=0.24\columnwidth,height=.21\linewidth]{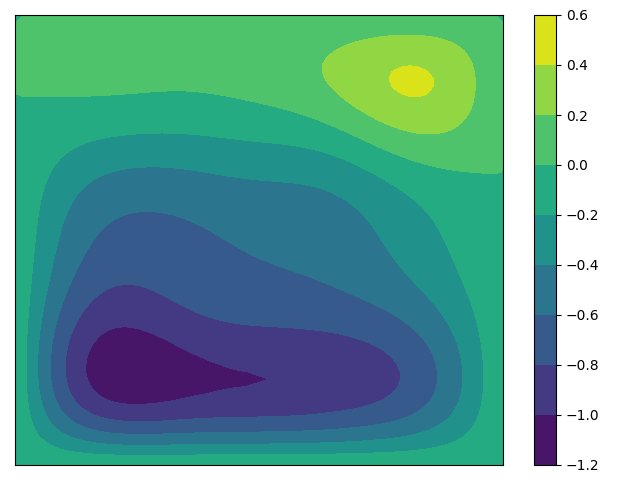}
        \includegraphics[width=.24\columnwidth,height=.21\linewidth]{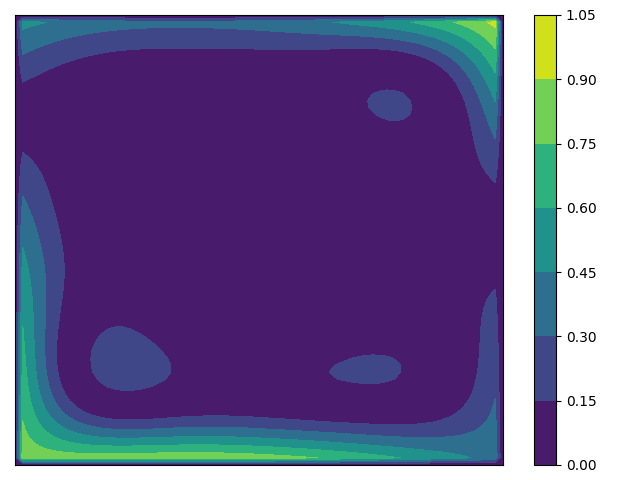}
    }
    \subfloat[The $e_f=\log{\frac{\|\phi_3-\phi_{3,rec}\|_2}{\|\phi_3\|_2}}$.]{
        \label{initial_value_distance_3}
        \includegraphics[width=0.24\columnwidth,height=.21\linewidth]{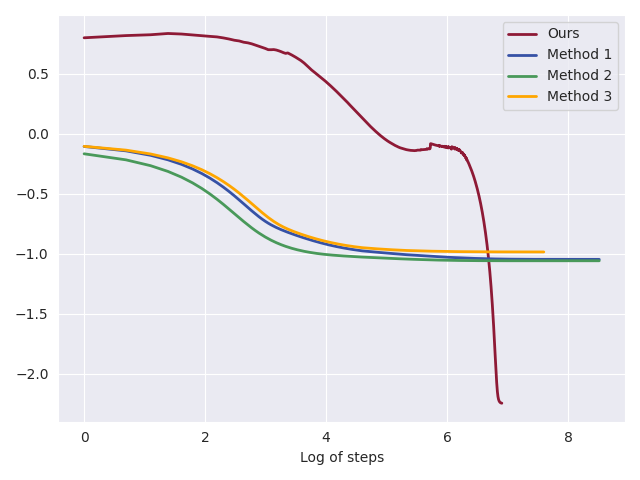}
    }
    \subfloat[The $e_u=\log{\frac{\|u^\delta(x,y,T;\phi_3)-u(x,y,T;\phi_{3,rec})\|_2}{\|u^\delta(x,y,T;\phi_3)\|_2}}$.]{
        \includegraphics[width=0.24\columnwidth,height=.21\linewidth]{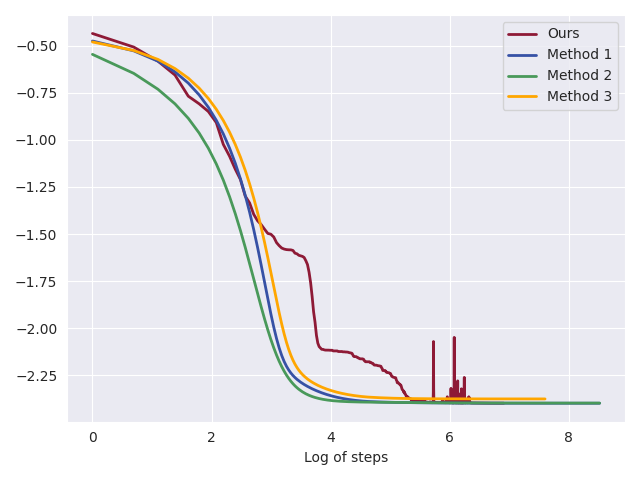}
    }
    \caption{The inversion results of initial value $\phi_3$. (a): the true $\phi_3$ and measurement data $u^\delta(x,y,T; \phi_3)$. (b)(c)(d)(e): the inversion $\phi_3$ and errors using different regularization methods. (f): the errors between the true $\phi_3$ and the inversion $\phi_{3,rec}$. (g): the errors between the measurement data $u^\delta(x,y,T; \phi_3)$ and $u(x,y,T; \phi_{3,rec})$. Here, we set the noise level $\varepsilon=0.05$ in (\ref{5.1}), the stopping parameter $\tau=1.01$ in discrepancy principle, the regularization parameter $\alpha=0.1$ in Tikhonov regularization.}
    \label{phi3_figure}
\end{figure}

 Similar to the previous inverse heat source problem, we invert the following three initial values from the measurement data $u^{\delta}(x,y,T)$ to validate the effectiveness of our method in inverse initial value problems for the heat equation:
\begin{itemize}
    \item \textbf{Case 1:} We select initial value $\phi_1$ share the same form of training data, but $\phi_1(x,y)$ is not included in the training set and not equal to $f_1$. The ground truth initial value $\phi_1$ is shown in Figure \ref{phi1_figure}.
    \item \textbf{Case 2:} We select $\phi_2(x,y)=4(x-x^2)\cos\pi y$, which is entirely different from the training data. The true initial value $\phi_2$ is depicted in Figure \ref{phi2_figure}.
    \item \textbf{Case 3:} We take $\phi_3(x,y)=\frac{7}{5}(x^2+(y-\frac{1}{3})^2)-1$, and the actual initial value $\phi_3$ is depicted in Figure \ref{phi3_figure}.
\end{itemize}


\begin{table}[ht]
    \begin{center}
    \begin{tabular}{llllllr}
    \toprule[1pt]
     &\multicolumn{2}{c}{$\phi_1$} & \multicolumn{2}{c}{$\phi_2$} & \multicolumn{2}{c}{$\phi_3$} \\
    \cline{2-3} \cline{4-5} \cline{6-7}
    method      & steps & relative $l_2$ error  & steps &relative $l_2$ error & steps & relative $l_2$ error\\ \hline
    Method 1  &  138 &  39.8\% & 86 & 36.4\% & 96 & 37.7\% \\ 
    Method 2     & 86 & 38.4\% & 76 & 36.0\% & 63 & 36.3\%\\ 
    Method 3     & 2000 & 40.5\% & 2000 & 36.4\% & 2000 & 37.3\%\\
    ODE-DPS(ours) & 1000 & \textbf{10.2\%} & 1000 & \textbf{6.9\%} & 1000 & \textbf{10.6\%} \\ 
    \bottomrule[1pt]
    \end{tabular}
    \end{center}
    \caption{The relative $l_2$ error and iteration steps of inverse initial value problem using different inversion methods. Methods 1, 2, and 3 correspond to Landweber iteration, improved Landweber iteration, and Tikhonov regularization, respectively. we set the noise level $\varepsilon=0.05$ in (\ref{5.1}), the stopping parameter $\tau=1.01$ in discrepancy principle, the regularization parameter $\alpha=0.1$ in Tikhonov regularization.}
    \label{initial_value_table}
\end{table}

\textit{\textbf{Results and discussion:}} 
Table \ref{initial_value_table} presents the relative $l_2$ error for different initial values $\phi$ using various regularization methods. Visual comparisons can be observed in Figures \ref{phi1_figure}-\ref{phi3_figure}. Our conclusions are summarized as follows:
\begin{enumerate}
    \item Similar to the experimental findings in identifying the heat source problem, the ODE-DPS inversion algorithm produced the most accurate numerical inversion results for all three numerical examples. Moreover, the relative $l_2$ error obtained with the three traditional regularization methods, is higher compared to the relative $l_2$ error in the inverse source problem. This discrepancy arises from the exponential decay of initial values in the heat equation, indicating that perturbations in the initial values have minimal impact on the solution at time $T$. Consequently, the inverse problem of recovering initial values becomes more challenging, further emphasizing the effectiveness of our algorithm.
    \item From Figures \ref{initial_value_distance_1}, \ref{initial_value_distance_2} and \ref{initial_value_distance_3}, we observe that the relative $l_2$ error of the numerical inversion results generated by traditional regularization methods decrease rapidly initially and then tends to level off, while the relative $l_2$ error produced by ODE-DPS show a consistent decrease and the inversion solution significantly superior to that obtained by traditional regularization methods.
\end{enumerate}

\subsection{Inverse Problem of Wave Equation}
\subsubsection{Problem settings}

To further demonstrate the effectiveness of our proposed inversion algorithm, we consider a wave equation as follows:
\begin{equation}
    \left\{\begin{aligned} 
        u_{tt}(x,y,t) &= a\Delta u(x,y,t) + f(x,y), (x,y,t)\in\Omega\times(0,T),\\
        u(x,y,t) &= g(x,y,t), (x,y,t)\in \partial\Omega\times (0,T),\\
        u(x,y,0) &= \phi(x,y), (x,y)\in\bar\Omega,\\
        u_t(x,y,0) &= \psi(x,y), (x,y)\in\bar\Omega,\\
    \end{aligned}\right.
    \label{wave}
\end{equation}
where $\Omega = [0,1]\times [0,1]$, the coefficient $a$ is a positive constant. Similar to the experimental settings discussed earlier for inverse problems in the heat equation, we focus solely on an inverse source problem utilizing measurement data $u^{\delta}(x,y,T)$ in this part to demonstrate the effectiveness of our method in wave equation inverse problems. 

\subsubsection{Implementation details}
To solve the wave equation (\ref{wave}), we employ a five-point explicit difference scheme for the forward problem simulation. We fix $T=1$ and set the time grid step size as $\Delta t=\frac{1}{100}$, and the spatial grid step size as $\Delta x=\Delta y=\frac{1}{63}$. For the noisy measurement data, we set $\varepsilon=0.05$. The functions for problem (\ref{wave}) are set as follows: the positive coefficient $a=0.05$, the boundary condition $g(x,y,t)=0$ for $(x,y,t)\in \partial\Omega\times(0,T)$, and the initial values $\phi(x,y)=0$ and $\psi(x,y)=0$ for $(x,y)\in\Omega$. The parameters of our ODE-DPS algorithm are set as $\gamma=0.65$, $\zeta=1.1$, $c=1e-5$, and $N=1000$. 
It is important to emphasize that, in this inverse source problem for the wave equation, there is no need to retrain the neural network. We can utilize the neural network $s_{\theta^*}(\mathbf{x}(t),t)$ previously trained for the inverse problems of the heat equation in the ODE-DPS inversion algorithm to identify the wave source problem.

\begin{figure}[htbp!]
    \centering
    \subfloat[The true source $f$(left) and measurements data(right).]{
        \includegraphics[width=0.24\columnwidth,height=.21\linewidth]{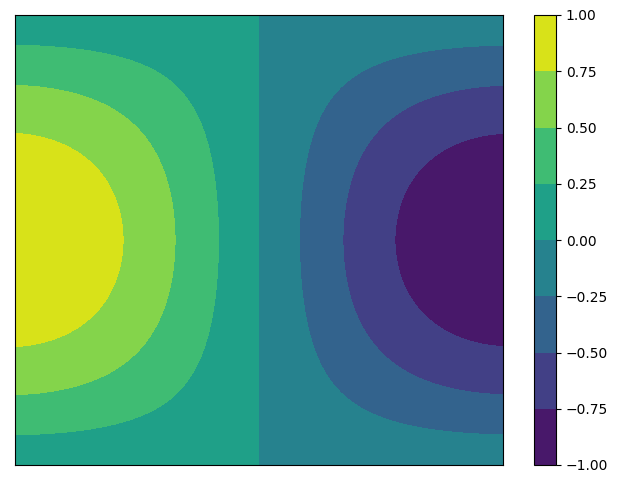}
        \includegraphics[width=.24\columnwidth,height=.21\linewidth]{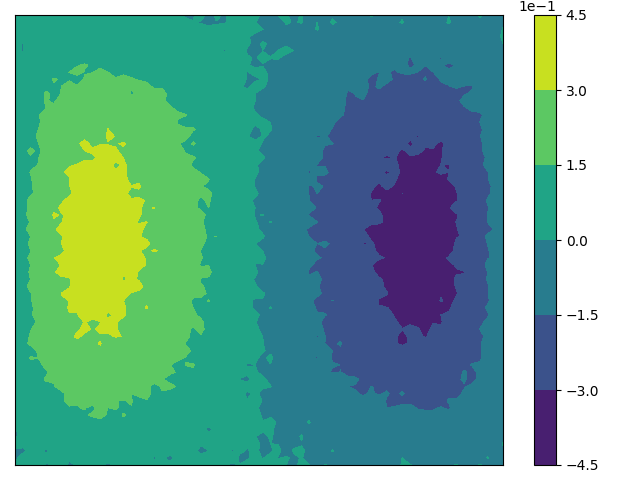}
    }
    \subfloat[The inversion $f$(left) and the error(right) by ODE-DPS.]{
        \includegraphics[width=0.24\columnwidth,height=.21\linewidth]{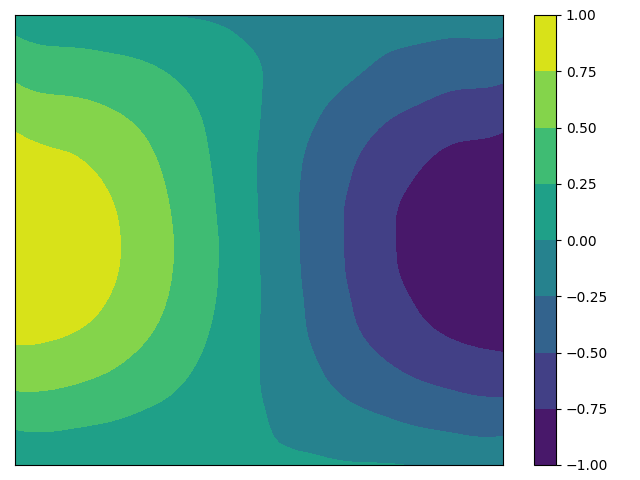}
        \includegraphics[width=.24\columnwidth,height=.21\linewidth]{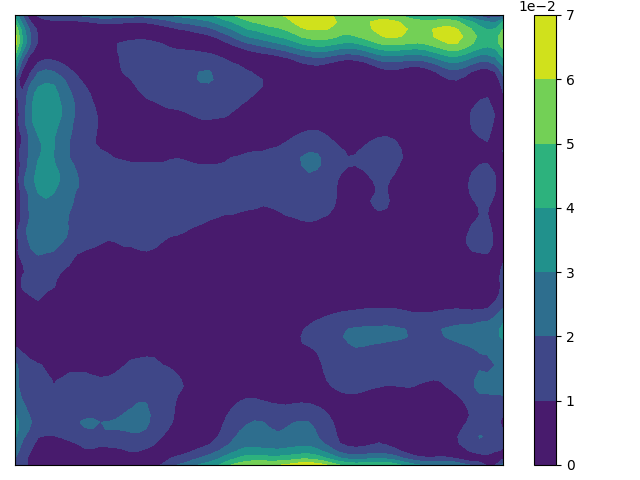}
    }
    
    \subfloat[The inversion $f$(left) and the error(right) by Method 1(Landweber iteration).]{
        \includegraphics[width=0.24\columnwidth,height=.21\linewidth]{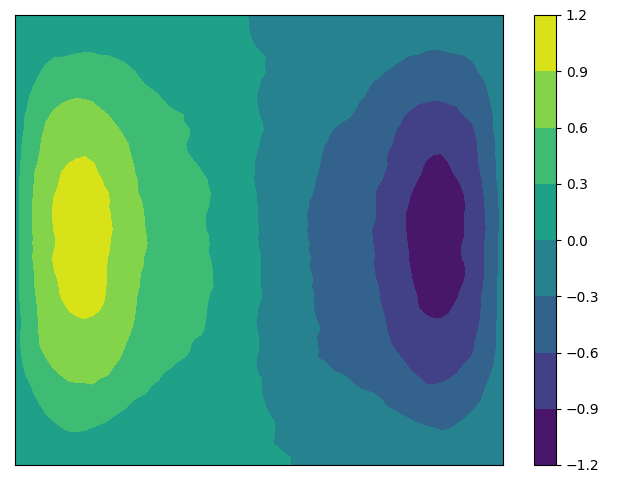}
        \includegraphics[width=.24\columnwidth,height=.21\linewidth]{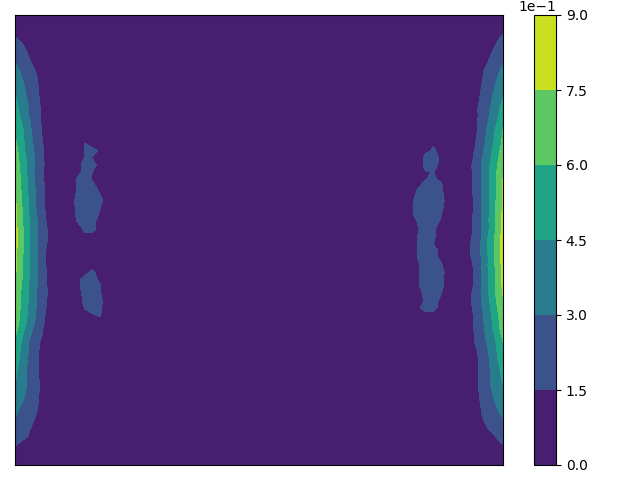}
    }
    \subfloat[The inversion $f$(left) and the error(right) by Method 2(Improved Landweber iteration).]{
        \includegraphics[width=0.24\columnwidth,height=.21\linewidth]{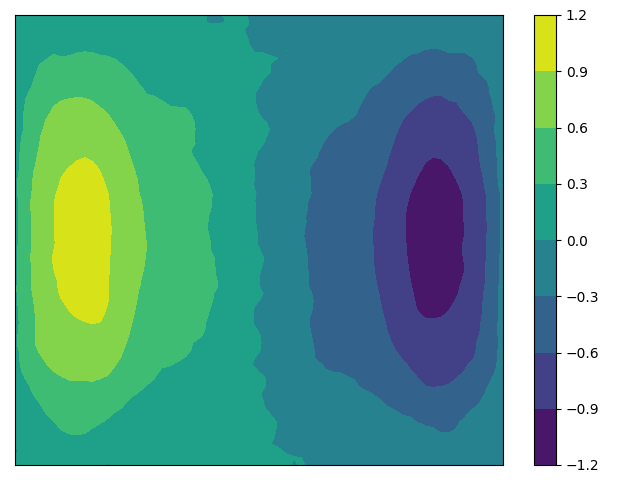}
        \includegraphics[width=.24\columnwidth,height=.21\linewidth]{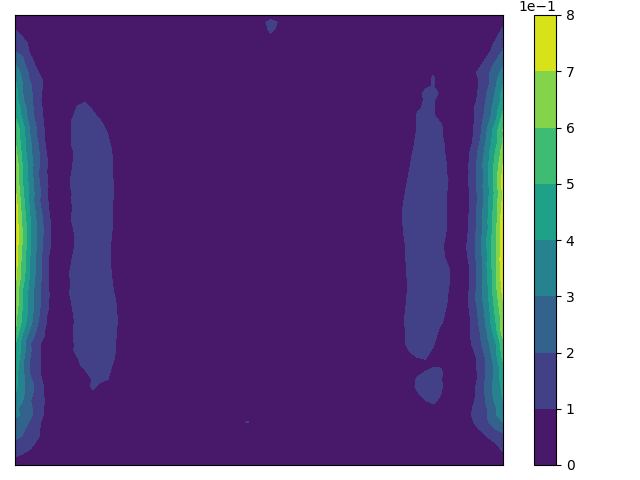}
    }

    \subfloat[The inversion $f$(left) and the error(right) by Method 3(Tikhonov regularization).]{
        \includegraphics[width=0.24\columnwidth,height=.21\linewidth]{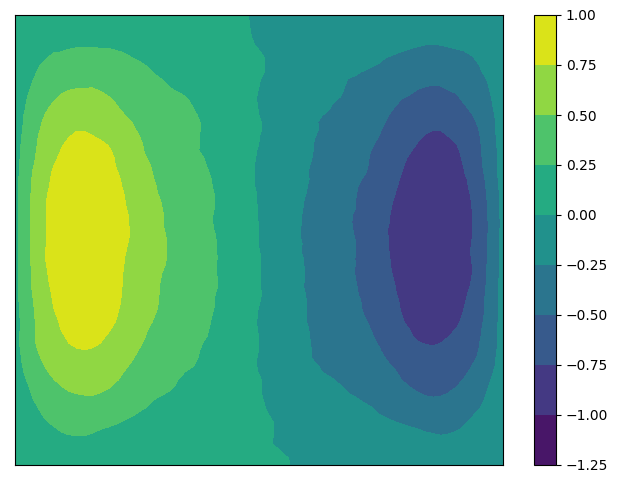}
        \includegraphics[width=.24\columnwidth,height=.21\linewidth]{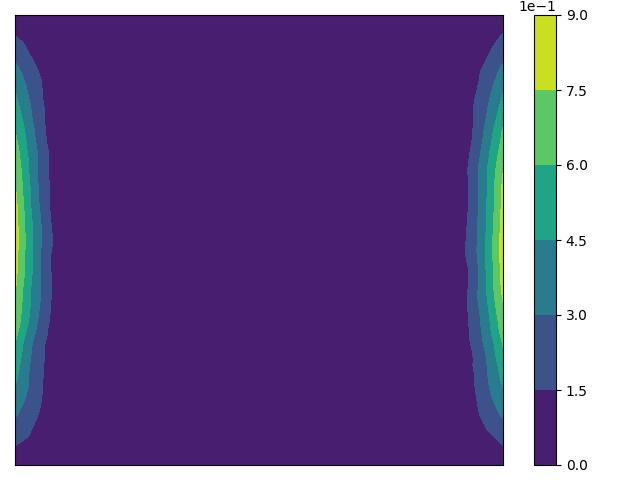}
    }
    \subfloat[The $e_f=\log\frac{\| f-f_{rec}\|_2}{\| f \|_2}$.]{\label{wave_fige}
        \label{wave_distance}
        \includegraphics[width=0.24\columnwidth,height=.21\linewidth]{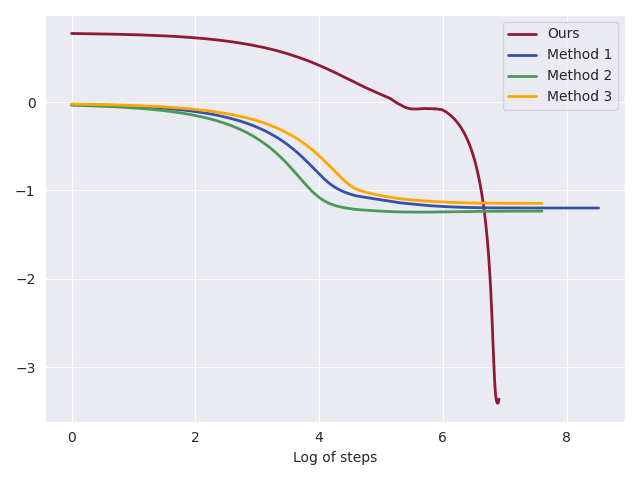}
    }
    \subfloat[The $e_u=\log{\frac{\| u^\delta(x,y,T;f)-u(x,y,T;f_{rec})\|_2}{\| u^\delta(x,y,T;f)\|_2}}$.]{\label{wave_figf}
        \includegraphics[width=0.24\columnwidth,height=.21\linewidth]{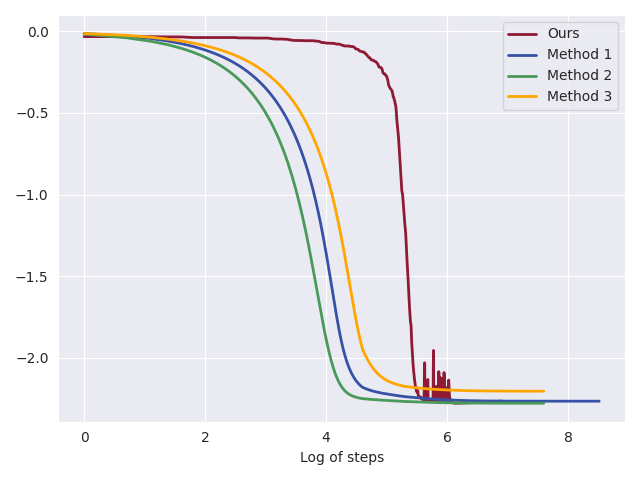}
    }
    \caption{The inversion results of source $f$ in wave equation. (a): the true $f$ and measurement data $u^\delta(x,y,T; f)$. (b)(c)(d)(e): the inversion $f$ and errors using different regularization methods. (f): the errors between the true $f$ and the inversion $f_{rec}$. (g): the errors between the measurement data $u^\delta(x,y,T;f)$ and $u(x,y,T;f_{rec})$. Here, we set the noise level $\varepsilon=0.05$ in (\ref{5.1}), the stopping parameter $\tau=1.01$ in discrepancy principle, the regularization parameter $\alpha=0.1$ in Tikhonov regularization.}
    \label{wave_figure}
\end{figure}

We choose $f(x,y)=\sin\pi x\cos\pi y$ as a representative to validate the effectiveness of our method in solving the inverse source problem of the wave equation. The ground truth source $f$ is shown in Figure \ref{wave_figure}. We apply our method on the inverse source problem in the wave equation and compare it with the traditional methods (Method 1, Method 2 and Method 3).

\textit{\textbf{Results and discussion}}:
Table \ref{wave_source_table} demonstrates the relative $l_2$ error for the inversion source $f$ using different regularization methods. The visual comparision can be found in Figure \ref{wave_figure}. We have the following findings:
\begin{enumerate}
    \item Similar to the experiments conducted for the heat equation, 
    all different regularization methods are capable of producing good numerical results. However, our proposed inversion algorithm outperforms the three compared traditional regularization methods in terms of lower relative $l_2$ errors in Table \ref{wave_source_table} and visually more accurate numerical reconstruction in Figure \ref{wave_figure}. However, our proposed inversion algorithm required more iterations compared to the first two methods, which mainly depends on the sampling step size of the scored-based diffusion model.
    \item From Figures \ref{wave_fige} and \ref{wave_figf}, it can be observed that the experimental results of the wave equation are remarkably similar to those of the heat equation, except that the numerical inversion results generated by Landweber iterative regularization method (Method 1) converge to a stable solution rather than exhibit a trend of decreasing relative $l_2$ error initially followed by an increase as the iteration progresses. However, this stable solution is obviously inferior to the solution obtained by our method.
\end{enumerate}

\begin{table}[htbp!]

    \begin{center}
    \begin{tabular}{llr}
    \toprule[1pt]
    method      & steps & relative $l_2$ error  \\
    \hline
    Method 1       &468      &30.5\%  \\ 
    Method 2       &142      &29.1\% \\ 
    Method 3      & 2000 & 31.8\%\\
    ODE-DPS(ours)  &1000     &\textbf{3.5\%}  \\
    \bottomrule[1pt]
    \end{tabular}
    \end{center}
    \caption{The relative $l_2$ error and iteration steps of inverse source problem of wave equation (\ref{wave}) using different inversion methods. Methods 1, 2, and 3 correspond to Landweber iteration, improved Landweber iteration, and Tikhonov regularization, respectively. We fix the maximum number of iterations for Methods 1, 2, and 3 at 2000. We set the noise level $\varepsilon=0.05$ in (\ref{5.1}), the stopping parameter $\tau=1.01$ in discrepancy principle, the regularization parameter $\alpha=0.1$ in Tikhonov regularization.}
    \label{wave_source_table}
\end{table}

\subsection{Algorithm Investigations}

\subsubsection{Hyperparameters Analysis}
In this part, we will discuss the impact of selecting different hyperparameters on the experimental results. We take the inverse source problem of the wave equation mentioned above as an example, by varying the step sizes $\zeta$ and its decay rate $\gamma$ in gradient descent step of our method.

From Table \ref{tbl:robustness1}, we can observe that within the range of gradient descent step sizes fluctuating between 0.7 to 1.5, and decay rates fluctuating between 0.45 to 0.85, ODE-DPS consistently achieves low-error solutions. The difference between the best and worst solutions is only around 3\%. Moreover, focusing on the range of decay rates $\gamma$ fluctuating between 0.55 to 0.85, we can see that the relative $l_2$ errors obtained by the ODE-DPS method are all below 4.7\%, with fluctuations in relative $l_2$ error not exceeding 1.1\%. This indicates that the algorithm's results provided by ODE-DPS are not strongly dependent on the choice of hyperparameters and can produce good results within an appropriate range of values.

\begin{table}[htbp!]
    \centering
    \setlength{\tabcolsep}{6mm}{
    \begin{tabular}{c| c c c c c }
    \toprule[1pt]
    \diagbox{$\zeta$}{$\gamma$} & 0.45 & 0.55 & 0.65 & 0.75 & 0.85 \\
    \hline
    0.7 & 6.6\% & 4.7\% & 3.6\% & 3.4\% & 3.7\% \\
    0.9 & 6.1\% & 4.3\% & 3.5\% & 3.5\% & 4.0\% \\
    1.1 & 5.7\% & 4.1\% & 3.5\% & 3.6\% & 4.3\% \\
    1.3 & 5.5\% & 4.0\% & 3.4\% & 3.8\% & 4.5\% \\
    1.5 & 5.2\% & 3.9\% & 3.6\% & 3.8\% & 4.6\% \\   
    \bottomrule[1pt]
    \end{tabular}}
    \caption{The relative $l_2$ error of the inversion source $f$ for wave equation (\ref{wave}) under different parameters. In algorithm \ref{al2}, we set $N=1000$, $c=1e-5$. For noisy measurement, we take $\varepsilon=0.05$ in (\ref{5.1}).}
    \label{tbl:robustness1}
\end{table}
\subsubsection{Robustness to Noise}

Next, we assess the algorithm's performance under varying levels of noise. We continue to use the inverse source problem of the wave equation mentioned above as an example, setting $\varepsilon$ in (\ref{5.1}) to $[0, 0.001, 0.005, 0.01, 0.05, 0.1]$. For all comparison methods, we fix the maximum number of iterations at 1000 and record the relative $l_2$ error under different noise levels, as presented in Table \ref{tbl:robustness3}. Our observations reveal that our method consistently outperforms the other three methods across all noise levels. Moreover, as the noise level increases, the error of the solutions generated by our method remains relatively stable, whereas the errors of the solutions obtained using the other three methods are significantly impacted. The data in Table \ref{tbl:robustness3} underscores the remarkable robustness of our proposed method to noise, surpassing all three traditional methods in performance.
\begin{table}[htbp!]
    \centering
    \setlength{\tabcolsep}{3mm}{
    \begin{tabular}{c| c c c c c c}
    \toprule[1pt]
    \diagbox{Methods}{$\varepsilon$} & 0 & 0.001 & 0.005 & 0.01 & 0.05 & 0.1\\
    \hline
    ODE-DPS(ours)  & \textbf{2.0\%} & \textbf{2.2\%} & \textbf{2.1\%} & \textbf{2.3\%} & \textbf{3.5\%}& \textbf{4.6\%}\\
    Method 1  & 25.2\% & 25.4\% & 26.8\% & 27.8\%& 30.5\%& 33.9\% \\
    Method 2  & 21.6\% & 21.8\% & 23.3\% & 24.7\% & 29.1\%& 30.9\%\\ 
    Method 3  & 27.0\% & 27.1\% & 28.0\% & 28.8\% & 31.8\%& 34.4\%\\ 
    \bottomrule[1pt]
    \end{tabular}}
    \caption{The relative $l_2$ error of the inversion source $f$ for wave equation (\ref{wave}) under different noise level. Methods 1, 2, and 3 correspond to Landweber iteration, improved Landweber iteration, and Tikhonov regularization, respectively. We fix the maximum number of iterations for all algorithms to 1000. We set the noise level $\varepsilon=0.05$ in (\ref{5.1}), the stopping parameter $\tau=1.01$ in discrepancy principle, the regularization parameter $\alpha=0.1$ in Tikhonov regularization.}
    \label{tbl:robustness3}
\end{table}

\subsection{Advantages and Limitations}
\textbf{Advantages:} In the experimental section, we employ the ODE-DPS inversion algorithm to solve three PDE inverse problems: the inverse source problem and inverse initial value problem in the heat equation, as well as the inverse source problem in the wave equation. We compare our proposed inversion method with traditional approaches in these experiments. Through extensive numerical experiments, we summarize the advantages of our algorithm as follows:
\begin{enumerate}
    \item Compared to traditional inversion methods, our proposed ODE-DPS inversion algorithm significantly enhances the inversion accuracy. These facts are evident from Tables \ref{heat_source_table}, \ref{initial_value_table}, \ref{wave_source_table}. For instance, in the inverse heat source problem, the lowest relative $l_2$ error obtained by the traditional regularization method across three examples is 24.1\%, while our ODE-DPS inversion algorithm achieves the highest relative $l_2$ error is 10.4\%. 
    \item Our inversion algorithm does not operate on a case-by-case area. That is, our proposed inversion method requires training the neural network $s_{\theta}(\bar{x},t)$ only once, and it can subsequently be applied to various inverse problem-solving tasks. This demonstrates that our deep learning-based inversion algorithm achieves nearly the same inversion efficiency as traditional inversion methods.
    \item Our inversion method demonstrates robustness to both the parameters within the algorithm and the level of noise present in the measurement data. This implies that, in experiments, we can achieve high-precision inversion results without the need for exhaustive parameter tuning. Moreover, the algorithm's robustness to measurement noise suggests that even in scenarios with significant noise in the measurement data, our proposed inversion method consistently produces accurate numerical results.
    \item Compared to the Landweber iteration regularization method, our ODE-DPS inversion method does not exhibit the property of semi-convergence. This fact is evident from the relative $l_2$ error plots of the experimental results.
\end{enumerate}

\textbf{Limitations:}
Due to the discretization step size dictated by the diffusion model's backward generation of the differential equation, our inversion method require $N=1000$ iterations to produce the final inversion outcome. This slightly increases the computational burden for solving some inverse problems. In the future, we aim to improve this aspect.

\section{Conclusion}\label{sec6}
In this paper, we introduce an ODE-based Diffusion Posterior Sampling (ODE-DPS) algorithm aimed at solving inverse problems associated with partial differential equations (PDEs). The method operates within a Bayesian inversion framework, transforming the problem of solving the posterior distribution into a process of generating samples from conditional distributions. Subsequently, it utilizes the score-based diffusion model for solving. To improve the inversion accuracy of this approach, we propose the ODE-DPS algorithm and furnish a detailed mathematical exposition thereof. Additionally, we conduct numerical experiments to illustrate the efficacy of our method. We anticipate that our proposed approach will be applied in various PDEs inverse problem-solving.

\bigskip
\noindent{\bf Acknowledgments}\\
Zheng Ma is supported by NSFC Grant No. 92270120 and No. 12201401. This work is partially supported by NSFC Grant No. 12031013 and the Strategic Priority Research Program of Chinese Academy of Sciences, XDA25010404. Additionally, it is also partially supported by Institute of Modern Analysis – A Shanghai Frontier Research Center.
\medskip

\clearpage
\medskip
\bibliographystyle{unsrt}
\bibliography{reference}

\clearpage
\appendix

\section{Training Algorithm}\label{A}

Below, we will provide some insights into the neural network training algorithm.
To initiate the noising process (\ref{2.2}), we fix the ending time $\hat{T}=1$ and set the time interval $\Delta t=\frac{1}{N}$. Then the time $\hat{T}$ is uniformly divided into $N$ parts, denoted as $t_i=i\Delta t, i=0,1,\cdots,N$. At each time point $t_i$, we determine the value of $\beta(t)$ and consider $\beta(t_i)\Delta t$ and $\mathbf{x}(t_i)$ as $\beta_i$ and $\mathbf{x}_i$, respectively. These values satisfy the condition $0<\beta_1<\beta_2<\cdots<\beta_N<1$. It's worth noting that $\beta_i$ can either be constant or learned through reparameterization. For the sake of simplicity, we assume $\beta_i$ to be fixed. Following this, for every data point $\mathbf{x}_0$ in the dataset, we add Gaussian noise based on the discretization of the stochastic differential equation (\ref{2.2}) as outlined below:
\begin{equation*}
\mathbf{x}_i - \mathbf{x}_{i-1} = -\frac{1}{2}\beta(t_i)\Delta t \mathbf{x}_{i-1}+\sqrt{\beta(t_i)\Delta t} z_i.
\end{equation*}
As interval $\Delta t\ll 1$, we can approximate $1-\frac{1}{2}\beta_i$ with $\sqrt{1-\beta_i}$ (Using Taylor's expansion) and get:
\begin{equation*}
\mathbf{x}_i \approx \sqrt{1-\beta_i}\mathbf{x}_{i-1}+\sqrt{\beta_i}z_i,
\end{equation*}
where $z_i\sim\mathcal{N}(0,I).$

From above discussion, we get the probability distributions of $\mathbf{x}_i$ given $\mathbf{x}_{i-1}$ as:
\begin{equation*}
p(\mathbf{x}_i|\mathbf{x}_{i-1}) = \mathcal{N}(\mathbf{x}_i;\sqrt{1-\beta_i}\mathbf{x}_{i-1}, \beta_iI).
\end{equation*}
Using the notation $\alpha_i:=1-\beta_i$ and $\bar{\alpha}_i:=\prod_{j=1}^i \alpha_i$ and the same derivations as described in reference \cite{ho2020denoising}, we have:
\begin{equation*}
p(\mathbf{x}_i|\mathbf{x}_0) = \mathcal{N}(\mathbf{x}_i;\sqrt{\bar{\alpha}_i}\mathbf{x}_0,(1-\bar{\alpha}_i)I).
\end{equation*}

As mentioned earlier, our aim is to train a neural network $s_\theta(\mathbf{x}(t), t)$, where $\theta$ represents the network parameters, to approximate 
$\nabla_{\mathbf{x}(t)}\log p(\mathbf{x}(t))$. This enables us to reverse the forward process and generate the data. Naturally, we have an optimization objective:
\begin{equation*}
J_1=\mathbb{E}_{t\sim U(0,1), \mathbf{x}(t)\sim p(\mathbf{x}(t))}[\| s_\theta(\mathbf{x}(t), t) - \nabla_{\mathbf{x}(t)}\log p(\mathbf{x}(t))\|_2^2].
\end{equation*}
However, the precise distribution $p(\mathbf{x}(t))$ is unknown and we have to construct an other objective equivalent to $J_1$ and could be calculated. We use denoising score matching \cite{vincent2011connection} to establish another objective as bellow:

$$J_2 = \mathbb{E}_{t\sim U(0,1),\mathbf{x}(t)\sim p(\mathbf{x}(t)|\mathbf{x}(0)), \mathbf{x}(0)\sim p(\mathbf{x}(0))}[\|s_\theta(\mathbf{x}(t), t)-\nabla_{\mathbf{x}(t)}\log p(\mathbf{x}(t)|\mathbf{x}(0))\|_2^2].$$
After discretization, $p(\mathbf{x}(t)|\mathbf{x}(0))$ can be written as $p(\mathbf{x}_i|\mathbf{x}_0),i=1,2,\cdots,N$. Besides, the equivalence of 
$J_1$ and $J_2$ is provided in the \ref{C}. 

Then, we can utilize optimization techniques to update parameter $\theta$ in order to minimize $J_2$ until convergence, thereby obtaining the optimal solution $\theta^{*}$ and the neural network $s_{\theta^*}(\mathbf{x}(t),t)$ to approximate the score $\nabla_{\mathbf{x}(t)} p(\mathbf{x}(t))$. We summarize the training process of the neural network $s_{\theta}(\mathbf{x}(t),t)$ in Algorithm \ref{al3}.

\begin{algorithm}[t]
    \textbf{Require}: $\mathbf{X}, N, \{\beta_i\}_{i=0}^N$, where $\mathbf{X}$ is a dataset; \\
    \textbf{Result}: $\theta$ \\
    $\alpha_i := 1-\beta_i,\bar{\alpha_i}:=\prod_{j=1}^{i}\alpha_i$\\
    \For{$\mathbf{x}_0\in \mathbf{X}$}{
    \For{$i=1$ \KwTo $N$}{
        $\mathbf{x}_i \sim \mathcal{N}(\mathbf{x}_i;\sqrt{\alpha_i}\mathbf{x}_{i-1}, \beta_iI)$\\
        $p(\mathbf{x}_i|\mathbf{x}_0)=\mathcal{N}(\mathbf{x}_i;\sqrt{\bar{\alpha}_i}\mathbf{x}_0, (1-\bar{\alpha}_i)I)$\\
    }
    \textbf{Repeat}\\
    $J_2=\mathbb{E}_{t\sim U(\epsilon, 1),\mathbf{x}(t)\sim p(\mathbf{x}(t)|\mathbf{x}(0)),\mathbf{x}(0)\sim p(\mathbf{x}(0))}[\|s_\theta(\mathbf{\mathbf{x}}(t),t)-\nabla_{\mathbf{x}(t)}\log p(\mathbf{x}(t)|\mathbf{x}(0))\|_2^2]$\\
    Update $\theta$ with Adam optimization algorithm.\\
    \textbf{Until} converged
    }
    \caption{Training algorithm.}
    \label{al3}
\end{algorithm}

\section{Proof of equivalence between $J_1$ and $J_2$.}
\label{C}
As shown in \ref{A}:
\begin{align*}
J_1&=\mathbb{E}_{t\sim U(0,1),~\mathbf{x}(t)\sim p(\mathbf{x}(t))}[\| s_\theta(\mathbf{x}(t), t) - \nabla_{\mathbf{x}(t)}\log p(\mathbf{x}(t))\|_2^2]\\ 
J_2 &= \mathbb{E}_{t\sim U(0,1),~\mathbf{x}(t)\sim p(\mathbf{x}(t)|\mathbf{x}(0)), ~\mathbf{x}(0)\sim p(\mathbf{x}(0))}[\|s_\theta(\mathbf{x}(t), t)-\nabla_{\mathbf{x}(t)}\log p(\mathbf{x}(t)|\mathbf{x}(0))\|_2^2].
\end{align*}
Since $t$ follows the same distribution in both $J_1$ and $J_2$, we only consider the equivalence of minimizing the following problem:
\begin{align*}
\hat{J}_1 &= \mathbb{E}_{\mathbf{x}(t)\sim p(\mathbf{x}(t))}[\| s_\theta(\mathbf{x}(t), t) - \nabla_{\mathbf{x}(t)}\log p(\mathbf{x}(t))\|_2^2]\\ 
\hat{J}_2 &= \mathbb{E}_{\mathbf{x}(t)\sim p(\mathbf{x}(t)|\mathbf{x}(0)), ~\mathbf{x}(0)\sim p(\mathbf{x}(0))}[\|s_\theta(\mathbf{x}(t), t)-\nabla_{\mathbf{x}(t)}\log p(\mathbf{x}(t)|\mathbf{x}(0))\|_2^2].
\end{align*}
To do it, we rewrite them as: 
\begin{align*}
\hat{J}_1 &= \mathbb{E}_{\mathbf{x}(t)\sim p(\mathbf{x}(t))}[\| s_\theta(\mathbf{x}(t))\|_2^2]-S_1(\theta)+C_1\\
\hat{J}_2 &= \mathbb{E}_{\mathbf{x}(t)\sim p(\mathbf{x}(t)|\mathbf{x}(0)), ~\mathbf{x}(0)\sim p(\mathbf{x}(0))}[\|s_\theta(\mathbf{x}(t))\|_2^2]-S_2(\theta)+C_2,
\end{align*}
where $C_1=\mathbb{E}_{\mathbf{x}(t)\sim p(\mathbf{x}(t))}[\|\nabla_{\mathbf{x}(t)}\log p(\mathbf{x}(t))\|_2^2]$ and $C_2=\mathbb{E}_{\mathbf{x}(t)\sim p(\mathbf{x}(t)|\mathbf{x}(0)),~ \mathbf{x}(0)\sim p(\mathbf{x}(0))}[\|\nabla_{\mathbf{x}(t)}\log p(\mathbf{x}(t)|\mathbf{x}(0))\|_2^2]$. Both $C_1$ and $C_2$ are constants that do not depend on the training parameter $\theta$. Besides, 
\begin{align*}
S_1(\theta)&=2\mathbb{E}_{\mathbf{x}(t)\sim p(\mathbf{x}(t))}[\left\langle s_\theta(\mathbf{x}(t), t),\nabla_{\mathbf{x}(t)}\log p(\mathbf{x}(t))\right\rangle]\\
S_2(\theta)&=2\mathbb{E}_{\mathbf{x}(t)\sim p(\mathbf{x}(t)|\mathbf{x}(0)), ~\mathbf{x}(0)\sim p(\mathbf{x}(0))}[\left\langle s_\theta(\mathbf{x}(t), t),\nabla_{\mathbf{x}(t)}\log p(\mathbf{x}(t)|\mathbf{x}(0))\right\rangle].
\end{align*}
We will show $S_1(\theta)=S_2(\theta)$ as bellow:
\begin{align*}
&\frac{1}{2}S_1(\theta)\\
&=\int_{\mathbf{x}(t)}p(\mathbf{x}(t))\left\langle s_\theta(\mathbf{x}(t),t), \nabla_{\mathbf{x}(t)} \log p(\mathbf{x}(t)) \right\rangle d\mathbf{x}(t)\\
&=\int_{\mathbf{x}(t)}\left\langle s_\theta(\mathbf{x}(t),t), \nabla_{\mathbf{x}(t)} p(\mathbf{x}(t))\right\rangle d\mathbf{x}(t)\\
&=\int_{\mathbf{x}(t)}\left\langle s_\theta(\mathbf{x}(t),t),  \int_{\mathbf{x}(0)}p(\mathbf{x}(0))\nabla_{\mathbf{x}(t)}p(\mathbf{x}(t)|\mathbf{x}(0))d\mathbf{x}(0) \right\rangle d\mathbf{x}(t)\\
&=\int_{\mathbf{x}(t)}\left\langle s_\theta(\mathbf{x}(t),t),  \int_{\mathbf{x}(0)}p(\mathbf{x}(0))p(\mathbf{x}(t)|\mathbf{x}(0))\nabla_{\mathbf{x}(t)}\log p(\mathbf{x}(t)|\mathbf{x}(0))d\mathbf{x}(0) \right\rangle d\mathbf{x}(t)\\
&=\int_{\mathbf{x}(t)}\int_{\mathbf{x}(0)}p(\mathbf{x}(0))p(\mathbf{x}(t)|\mathbf{x}(0))\left\langle s_\theta(\mathbf{x}(t),t),  \nabla_{\mathbf{x}(t)}\log p(\mathbf{x}(t)|\mathbf{x}(0)) \right\rangle d\mathbf{x}(0)d\mathbf{x}(t)\\
&=\mathbb{E}_{\mathbf{x}(t)\sim p(\mathbf{x}(t)|\mathbf{x}(0)),~ \mathbf{x}(0)\sim p(\mathbf{x}(0))}[\left\langle s_\theta(\mathbf{x}(t), t),\nabla_{\mathbf{x}(t)}\log p(\mathbf{x}(t)|\mathbf{x}(0))\right\rangle]\\
&=\frac{1}{2}S_2(\theta).
\end{align*}
Besides, we note that $$\mathbb{E}_{\mathbf{x}(t)\sim p(\mathbf{x}(t))}[\| s_\theta(\mathbf{x}(t))\|_2^2]=\mathbb{E}_{\mathbf{x}(t)\sim p(\mathbf{x}(t)|\mathbf{x}(0)), ~\mathbf{x}(0)\sim p(\mathbf{x}(0))}[\|s_\theta(\mathbf{x}(t))\|_2^2].$$
Based on the above discussion, we have $J_2=J_1-C_1+C_2$.
We have thus shown that the two optimization objectives $J_1$ and $J_2$ are
equivalent.
\clearpage
\section{The Detailed Structure of the Neural Network}\label{B}
\begin{center}
\includegraphics[width=1.4\columnwidth,height=0.6\linewidth,angle=90]{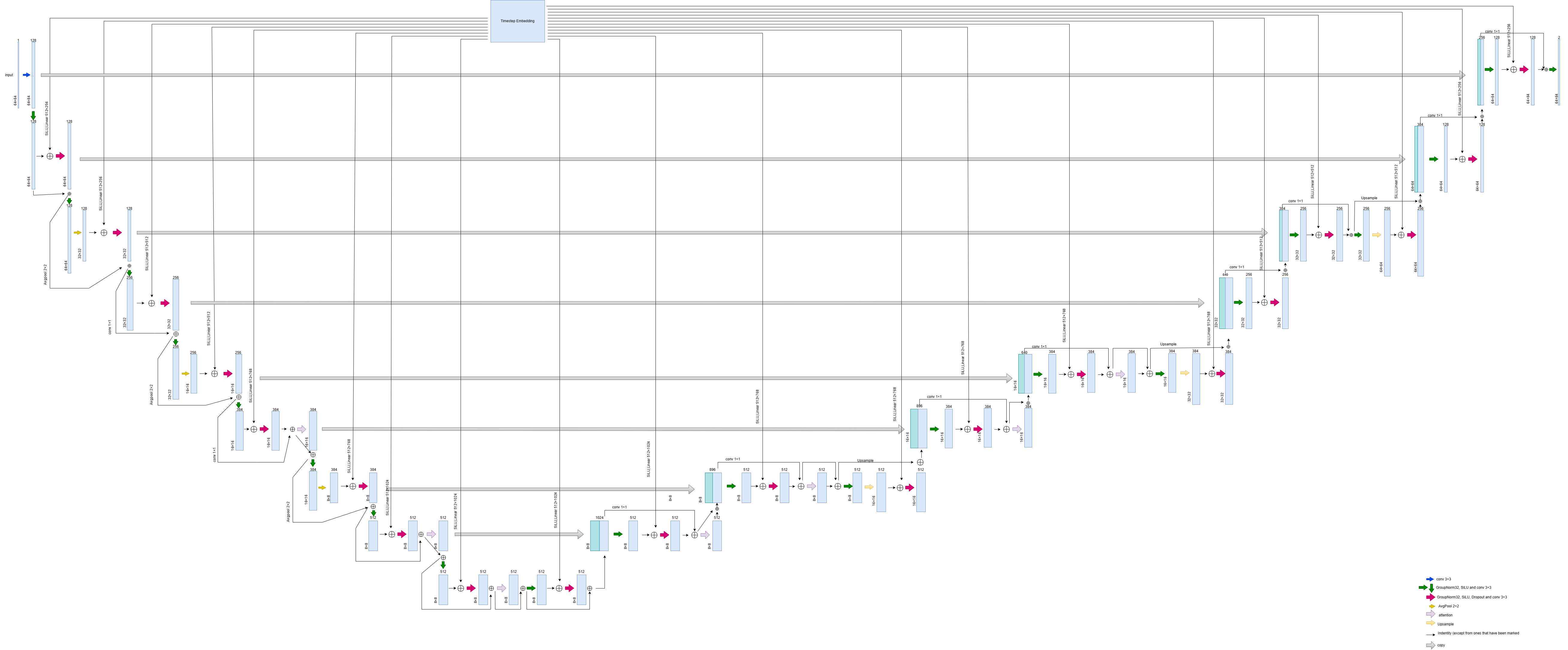}
\end{center}

\end{document}